\documentclass[11pt,letterpaper]{amsart}
 \usepackage{stmaryrd} 
 \usepackage{mathrsfs}
 \usepackage{amsmath}
 \usepackage{amsthm}
 \usepackage{amsfonts}
 \usepackage{amssymb}
 \usepackage{amsxtra}
 \usepackage[dvips]{graphicx}
 \usepackage{wrapfig}
 \usepackage{layout}
 \usepackage{alltt}
 \usepackage{xypic}
 \input xy
 \xyoption{all}

\DeclareMathAlphabet{\mathpzc}{OT1}{pzc}{m}{it}

\usepackage[pdftex,
                paper=a4paper,
                portrait=true,
                textwidth=450pt,
                textheight=673pt,
                tmargin=3cm,
                marginratio=1:1]{geometry}

\def\Cx{{\mathbb C}}

\def\Nx{{\mathbb N}}
\def\Px{{\mathbb P}}
 
\def\Rx{{\mathbb R}}

\def\Zx{{\mathbb Z}}

\def\id{{\rm id}}

\newcommand{\As}{\mathscr{A}}

\newcommand{\Fs}{\mathscr{F}}

\newcommand{\Hs}{\mathscr{H}}

\newcommand{\Ls}{\mathscr{L}}
\newcommand{\Ms}{\mathscr{M}}

\newcommand{\Os}{\mathscr{O}}
\newcommand{\Ps}{\mathscr{P}}

\newcommand{\Ss}{\mathscr{S}}

\newcommand{\Us}{\mathscr{U}}
\newcommand{\Vs}{\mathscr{V}}

\newcommand{\Xs}{\mathscr{X}}
\newcommand{\Ys}{\mathscr{Y}}

\newcommand{\Ld}{\mathbf{L}}
\newcommand{\Rd}{\mathbf{R}}
\newcommand{\Dd}{\mathbf{D}}

\def\hfrak{\mathfrak{h}}

\def\ufrak{\mathfrak{u}}

\newcommand{\fl}{\rightarrow}
\newcommand{\lfl}{\longrightarrow}

\newcommand\umod{\operatorname{mod \textendash \!}}

\newcommand\uMat{\operatorname{Mat}}
\newcommand\uEnd{\operatorname{End}}

\newcommand\uCoh{\operatorname{Coh}}

\newcommand\uCrit{\operatorname{Crit}}

\newcommand\uHom{\operatorname{Hom}}

\newcommand\uProj{\operatorname{Proj}}
\newcommand\urank{\operatorname{rank}}

\newcommand\uSpec{\operatorname{Spec}}

\newcommand\uRep{\operatorname{Rep}}

\newcommand\ugr{\operatorname{gr}}
\newcommand\ugldim{\operatorname{gldim}}

\newcommand\udeg{\operatorname{deg}}
\newcommand\utr{\operatorname{tr}}

\newcommand\uGL{\operatorname{GL}}

\newcommand\uSL{\operatorname{SL}}

\newcommand\uU{\operatorname{U}}

\newcommand\uch{\operatorname{ch}}

\newcommand\udimvec{\operatorname{\underline{dim}}}
\newcommand\usm{\operatorname{\#}}

\theoremstyle{plain}
 \newtheorem{thm}{Theorem}[section]
 \newtheorem{lem}[thm]{Lemma}
 \newtheorem{prop}[thm]{Proposition}
 \newtheorem{cor}[thm]{Corollary}
 \theoremstyle{definition}
 
 \newtheorem{exam}[thm]{Example}

\newtheorem*{ack}{{\footnotesize Acknowledgments}}
 \theoremstyle{remark}

\author{Alexander Quintero V\'{e}lez and Alex Boer}
\email{quintero@math.uu.nl, boer@math.uu.nl}
\address{Mathematisch Instituut\\
Universiteit Utrecht\\
P.O.~Box 80.010, NL-3508 TA Utrecht\\
Nederland}

\title{Noncommutative resolutions of ADE fibered Calabi-Yau threefolds}
\begin{document}
\subjclass[2000]{16G20, 14A22}
\keywords{Representation of quivers, moduli spaces, noncommutative resolutions}
\begin{abstract} 
In this paper we construct noncommutative resolutions of a certain class of Calabi-Yau threefolds studied in \cite{CKV01} by F.~Cachazo, S.~Katz and C.~Vafa. The threefolds under consideration are fibered over a complex plane with the fibers being deformed Kleinian singularities. The construction is in terms of a noncommutative algebra introduced by V.~Ginzburg in \cite{Gin06}, which we call the ``$N=1$ ADE quiver algebra''.  
\end{abstract}

\maketitle
\section*{Introduction}
In recent years, there has been a great deal of interest in noncommutative algebra in connection with algebraic geometry, particularly in the study of singularities and their resolutions. The underlying idea in this context is that the resolutions of a singularity are closely linked to the structure of a noncommutative algebra.

The case of Kleinian singularities $X=\Cx^2/G$, for $G$ a finite subgroup of $\uSL(2,\Cx)$, was the first non-trivial example of this phenomenon, studied in \cite{CS98}. It was shown that the minimal resolution of $X$ is a moduli space of representations of the preprojective algebra associated to the action of $G$. This preprojective algebra is known to be Morita equivalent to the skew group algebra $\Cx[x,y] \usm G$, so we could alternatively use this algebra to construct the minimal resolution of $X$. Later,  M.~Kapranov and E.~Vasserot \cite{KV00} showed that there is a derived equivalence between $\Cx[x,y] \usm G$ and the minimal resolution of $X$. A similar statement was established by T.~Bridgeland, A.~King and M.~Reid \cite{BKR01} for crepant resolutions of quotient singularities $X=\Cx^3/G$ arising from a finite subgroup $G \subseteq \uSL(3,\Cx)$. In this case the crepant resolution of $X$ is realized as a moduli space of representations of the McKay quiver associated to the action of $G$, subject to a certain natural commutation relations (see \cite{CrIs04} and \S 4.4 of \cite{Gin06}). 

Various steps in the direction mentioned above have been taken in a series of papers \cite{BrSym05,Br06,GS04,Sz07,Wem07}, where several concrete examples have been discussed. More abstract approaches have also been put forward in \cite{Br99,BBS03}. The lesson to be drawn from these works is that for some singularities it is possible to find a noncommutative algebra $A$ such that the representation theory of this algebra dictates in every way the process of resolving these singularities. More precisely, it is shown that:
\begin{enumerate}
\item[{\tiny \textbullet}]the centre of $A$ is the coordinate ring of the singularity;
\item[{\tiny \textbullet}]the algebra $A$ is finitely generated as a module over its centre; 
\item[{\tiny \textbullet}]resolutions of the singularity are realized as moduli spaces of representations of $A$; 
\item[{\tiny \textbullet}]the category of finitely generated modules over $A$ is derived equivalent to the category of coherent sheaves on an appropriate resolution.
\end{enumerate}
Following the terminology of M.~Van den Bergh (cf.~\cite{Bergh02,Bergh04}) we may think of $A$ as a ``noncommutative resolution''.

This phenomenon also appears naturally in string theory in the context of ``geometric engineering''. There the singularity $X$ should be a Calabi-Yau threefold and one studies Type IIB string theory compactified on $X$. It turns out that a collection of D-branes located at the singularity gives rise to a noncommutative algebra $A$, which can be described as the path algebra of a quiver with relations. For a fixed quiver $Q$, this construction only depends on a ``noncommutative function'' called the superpotential.  As a consequence, the aforementioned derived equivalence establishes a correspondence between two different ways of describing a D-brane:~as an object of the derived category of coherent sheaves on a crepant resolution of $X$ and as a representation of the quiver $Q$.

Now, let us explain the situation on which we will focus. We shall study a kind of singular Calabi-Yau threefolds obtained by fibering the total space of the semi-universal deformation of a Kleinian singularity over a complex plane, subsequently termed ADE fibered Calabi-Yau threefolds. They have been defined and studied in the work of F.~Cachazo, S.~Katz and C.~Vafa \cite{CKV01} from the point of view of $N=1$ quiver gauge theories. The quiver diagrams of interest here are the extended Dynkin quivers of type $A$, $D$ or $E$. Following \cite{Gin06}, for such a quiver $Q$, we associate a noncommutative algebra $\mathfrak{A}^{\tau}(Q)$ which we call the ``$N=1$ ADE quiver algebra''. The choice of $\tau$ is encoded in the fibration data. The goal of this paper is to show that the $N=1$ ADE quiver algebra realizes a noncommutative resolution of the ADE fibered Calabi-Yau threefold associated with $Q$ and $\tau$. The proof of this result depends on two ingredients. On the one hand, we use the results in \cite{CS98}, on the construction of deformations of Kleinian singularities and their simultaneous resolutions in terms of $Q$. On the other hand, we use the results in \cite{CrH98} to construct a Morita equivalence between $\mathfrak{A}^{\tau}(Q)$ and a noncommutative crepant resolution $A^{\tau}$ in the sense of Van den Bergh; this allows us to use the techniques developed in \cite{Bergh02} to show that the derived category of finitely generated modules over $\mathfrak{A}^{\tau}(Q)$ is equivalent to the derived category of coherent sheaves on any small resolution of the ADE fibered Calabi-Yau threefold.

Some related results using different methods were obtained by B.~Szendr{\H{o}}i in \cite{Sz08}. He considers threefolds $X$ fibered over a general curve $C$ by ADE singularities and shows that D-branes on a small resolution of $X$ are classified by representations with relations of a Kronheimer-Nakajima-type quiver in the category $\uCoh(C)$ of coherent sheaves on $C$. The correspondence is given by a derived equivalence between a small resolution of $X$ and a sheaf of noncommutative algebras on $C$. In particular, there is a substantial overlap between Section~\ref{sec3} of this paper and the results of Ref.~\cite{Sz08}. 

The structure of the paper is as follows. Section~\ref{sec1} will be devoted to setting up the physical and mathematical context of our work. Even though many of the statements outlined in this section do not constitute rigorous mathematics, they provide motivation and background to what follows. In Section~\ref{sec2} we define the $N=1$ ADE quiver algebra $\mathfrak{A}^{\tau}(Q)$, describe some of its basic structure and prove that small resolutions of ADE fibered Calabi-Yau threefolds are realized as moduli spaces of representations of $\mathfrak{A}^{\tau}(Q)$. We conclude with a discussion of the derived equivalence between $\mathfrak{A}^{\tau}(Q)$ and a small resolution of an ADE fibered Calabi-Yau threefold in Section~\ref{sec3}.

\begin{ack}
{\footnotesize We would like to thank Bal\'{a}zs Szendr{\H{o}}i for helpful remarks and e-mail correspondence. A.Q.V. is grateful to Jan Stienstra for his constant guidance and also for some essential conversations. A.Q.V. is also grateful to Tom Bridgeland, Michel Van den Bergh and Marjory Jane Macleod for valuable discussions related to this work. A.B. thanks A.Q.V. for introducing him to this problem.}
\end{ack}

\section{Physical and mathematical context}\label{sec1}
This section is a digression providing a general context and motivation for what we are doing. The setup for our discussion is the reverse geometric engineering of singularities.

Let us start with some background about D-branes. Recall that type II superstrings are described by maps from a Riemann surface $\Sigma$, the ``worldsheet'' as it is called, to a ten-dimensional ``spacetime'' manifold $M$. In the simplest instance, a D-brane is a submanifold of $M$ on which open strings can end. This means that if a D-brane is present, then one needs to consider maps from a Riemann surface with boundaries to $M$ such that the boundaries are mapped to a given submanifold $W \subseteq M$. In this case one says that there is a D-brane wrapped on $W$. If $W$ is connected and has dimension $p+1$ we will refer to the brane as a D$p$-brane. $W$ itself is referred to as the ``worldvolume'' of the D-brane. In what follows, we assume that the underlying ten-dimensional space $M$ is decomposed as $M=\Rx^{1,3} \times X$ where $\Rx^{1,3}$ denotes the four-dimensional Minkowski space and $X$ is a six-dimensional space given by a Calabi-Yau threefold. We will specialize to D-branes in $M$ whose worldvolume is of the form $W=\Rx^{1,3} \times S$ with $S \subseteq X$. Forgetting about the manifold $\Rx^{1,3}$ for the moment, we will speak of D-branes wrapping $S$.

D-branes are generally more complicated objects than just submanifolds in an ambient spacetime, because in string theory they are realized as boundary conditions for a certain auxiliary quantum field theory on the Riemann surface $\Sigma$. More concretely, the data specifying the boundary conditions for the auxiliary theory on $\Sigma$ include a choice of a rank $r$ vector bundle $E$ on $S$ and a connection on it. From a physical viewpoint, such bundle should be thought of as $r$ coincident D-branes wrapped on the same submanifold $S$. When these facts are properly taken into account, it turns out that the dynamics on the worldvolume $S$ is an $N=1$ supersymmetric gauge theory with gauge group $\uU(r)$.

For the present discussion we will only consider D-branes in the open string topological B-model, in which case the submanifolds are complex submanifolds and the vector bundles are holomorphic. The above point of view can be generalised if one takes into account that a holomorphic vector bundle defined on a complex submanifold $S$ defines a coherent sheaf $i_* E$ (with $i\colon S \hookrightarrow X$ being the inclusion map). We are thus lead to consider coherent sheaves with support on a submanifold of $X$. However, more generally one would like to describe collections of D-branes and anti-D-branes.\footnote{An anti-D-brane has all the same physical properties as an arbitrary D-brane, modulo the fact that they try to annihilate each other.} As explained in \cite{Asp05}, among many other references, this forces us to consider not only coherent sheaves but complexes of coherent sheaves. Furthermore, maps between complexes are represented by tachyons and localization on quasi-isomorphisms is expected to be realized by renormalization group flow. Assembling this information, the picture that emerges is that D-branes do correspond to objects in the derived category of coherent sheaves on $X$.

Now let us briefly discuss the geometric engineering of gauge theories. The idea behind geometric engineering is to look at the gauge theories that arise on D-branes at singularities. To be more concrete, we take spacetime to be $M=\Rx^{1,3} \times X$ where $X$ is a Calabi-Yau threefold with an isolated singularity at $P \in X$ and consider a D$3$-brane wrapped on $\Rx^{1,3} \times \{ P \}$. In terms of the derived category this D$3$-brane is represented by the skyscraper sheaf $\Os_P$ on $X$. We want to determine the gauge theory on the worldvolume of such D$3$-brane. The crucial thing to note here is that $\Os_P$ is marginally stable to decay into a collection of so-called ``fractional'' branes $E_i$.\footnote{The existence of this decay is argued for in \cite{Asp06}. In modern language, this construction amounts to introducing a notion of Bridgeland stability for D-branes; see, for example, \cite{Bergman08}.} Each fractional brane may appear with multiplicity\footnote{The multiplicities $\alpha_i$ are uniquely determined via the condition $\sum_i \alpha_i \uch(E_i)=\uch(\Os_P)$ (as asserted in Eq.~(2.1) of \cite{Wijn08}).} $\alpha_i$ and so is associated to a factor of $\uU(\alpha_i)$ in the worldvolume gauge theory. Most notably, these branes are bound together by open string excitations corresponding to classes of morphisms $\phi_{ij} \in \uHom(E_i,E_j)$. As a result, the effective dynamics on the worldvolume $\Rx^{1,3}$ of the D$3$-brane is an $N=1$ supersymmetric gauge theory whose matter content can be conveniently encoded in a quiver with relations coming from a ``superpotential''. This will therefore be what is known as an $N=1$ quiver gauge theory.  

We can be somewhat more precise about this. Let $Q$ be a quiver with vertex set $I$ and denote by $\Cx Q$ the corresponding path algebra. A {\it superpotential} is a formal sum of oriented cycles on the quiver, i.e.~an element of the vector space $\Cx Q / [\Cx Q,\Cx Q]$. On this space we can define for every arrow $a$ a ``derivation'' $\partial_a$ that takes any occurrences of the arrow in an oriented cycle and removes them leading to a path from the head of $a$ to its tail. An $N=1$ {\it quiver gauge theory} consists of a quiver $Q$ together with a choice of a superpotential $W$ and a dimension vector $\alpha \in \Nx^{I}$. Given an $N=1$ quiver gauge theory, one can construct the $N=1$ {\it quiver algebra}
$$
A=\Cx Q / \left(\partial_a W \mid a \in Q \right).
$$ 

We now explain how this data gives rise to an $N=1$ supersymmetric gauge theory on $\Rx^{1,3}$. Suppose we are given an $N=1$ quiver gauge theory $(Q,W,\alpha)$. Representations of the quiver $Q$ of dimension vector $\alpha$ are given by elements of the vector space
$$
\uRep(Q,\alpha)=\bigoplus_{a \in Q} \uHom_{\Cx}(\Cx^{\alpha_{t(a)}},\Cx^{\alpha_{h(a)}})
$$
where $h(a)$ and $t(a)$ denote the head and tail vertices of an arrow $a$. The isomorphism classes correspond to orbits of the group $\uGL(\alpha)=\prod_{i \in I}\uGL(\alpha_i,\Cx)$ acting by conjugation. We define a Hermitian inner product on $\uRep(Q,\alpha)$  via the trace form $(x,y)=\sum_{a \in Q} \utr (x_a y^*_a)$, where $*$ denotes the adjoint map. Let $\uU(\alpha)$ denote the product of unitary groups $\prod_{i \in I} \uU(\alpha_i)$. This is a maximal compact subgroup of $\uGL(\alpha)$ and acts on $\uRep(Q,\alpha)$ preserving the Hermitian structure. The corresponding moment map $\mu_{\alpha}\colon \uRep(Q,\alpha) \fl \bigoplus_{i \in I}\ufrak(\alpha_i)$ is given by 
$$
\mu_{\alpha}(x)_i=\sqrt{-1}\Bigg(\sum_{h(a)=i} x_a x_{a}^* -  \sum_{t(a)=i}
x_{a}^*x_a\Bigg).
$$
Now consider the superpotential $W$. Recall that it is required to be a sum of oriented cycles in $Q$.  If $W=\sum_{i_1, \dots, i_{r}}a_{i_1} \cdots a_{i_r}$, then the function $W_{\alpha}\colon \uRep(Q,\alpha) \fl \Cx$ given by
$$
W_{\alpha}(x)=\sum_{i_1, \dots, i_{r}}\utr(x_{a_{i_1}} \cdots x_{a_{i_r}})
$$
is invariant under the action of $\uGL(\alpha)$, and thus also $\uU(\alpha)$. The upshot of all this is that an $N=1$ quiver gauge theory gives rise to a quadruple $(\uU(\alpha),\uRep(Q,\alpha),\mu_{\alpha},W_{\alpha})$. According to \cite[Supersolutions, \S 6.2]{DelFreed99}, the latter is the data one needs to specify an $N=1$ supersymmetric gauge theory on $\Rx^{1,3}$. The gauge group is given by $\uU(\alpha)$. The field content of the theory associated with the quiver is encoded as follows. We associate to each vertex $i$ an $N=1$ vector multiplet $\As_i=(A_i,\lambda_i,D_i)$ and to each arrow $a$ an $N=1$ chiral multiplet $\Phi_a=(\phi_a,\psi_a,F_a)$. In the first (vector) multiplet, $A_i$ is a connection on some principal $\uU(\alpha_i)$ bundle over $\Rx^{1,3}$, $\lambda_i$ is spinor with values in the adjoint bundle and $D_i$ is an auxiliary field with values in the adjoint bundle. Letting $R_i$ be the fundamental $\alpha_i$-dimensional representation of $\uU(\alpha_i)$, there is an associated vector bundle $E_i$. If $R_i^*$ is the dual representation of $R_i$ \textendash sometimes called ``antifundamental representation''\textendash~then the associated vector bundle is the dual bundle to $E_i$. In the chiral multiplet, $\phi_a$ is a section of $E_{t(a)}^* \otimes E_{h(a)}$, $\psi_a$ is a spinor with values in $E_{t(a)}^* \otimes E_{h(a)}$, and $F_a$ is an auxiliary field with values in $E_{t(a)}^* \otimes E_{h(a)}$. The chiral multiplets are therefore often called ``bifundamental fields''. The explicit Lagrangian, together with the relevant supersymmetry transformations, is given in \cite[Theorem~6.33]{DelFreed99}.  

Now we give a description of the moduli space of classical vacua of an $N=1$ quiver gauge theory. Stealing a look at Theorem~6.33 of \cite{DelFreed99}, the moduli space of classical vacua of such a quiver gauge theory is
$$
\Ms_{0}=\mu_{\alpha}^{-1}(0) \cap \uCrit(W_{\alpha}) / \uU(\alpha).
$$ 
Here `$\uCrit$' denotes the set of critical points. In general, $W_{\alpha}$ drops to a regular function $\overline{W}_{\alpha}$ on the symplectic quotient $\mu_{\alpha}^{-1}(0) / \uU(\alpha)$ and the moduli space of vacua is the set of critical points of $\overline{W}_{\alpha}$ on this quotient. There is another mathematical interpretation of this process, as a quotient in the sense of GIT: we complexify the group $\uU(\alpha)$ to $\uGL(\alpha)$, and consider the action of $\uGL(\alpha)$ on $\uRep(Q,\alpha)$. It turns out that we can identify $\Ms_{0}$ with the set of critical points of $W_{\alpha}$ on the affine quotient variety $\uRep(Q,\alpha) \sslash \uGL(\alpha)$. 

Let us restate the above in terms of the $N=1$ quiver algebra $A$. Let $\uRep(A,\alpha)$ denote the closed subspace of $\uRep(Q,\alpha)$ corresponding to representations for $A$. The group $\uGL(\alpha)$ acts naturally on this variety, and the orbits correspond to isomorphism classes of representations. It is pointed out in \cite[Proposition~3.8]{Seg08} that the set of critical points of $W_{\alpha}$ is precisely $\uRep(A,\alpha)$. Altogether this implies that the moduli space of vacua $\Ms_{0}$ admits an alternate presentation as an affine quotient $\uRep(A,\alpha) \sslash \uGL(\alpha)$. At least in passing, we should mention that this construction has a noncommutative-geometric interpretation in which $A$ is viewed as the noncommutative coordinate ring of the ``critical locus'' of $W$; see, for instance, \cite{CEG07}.  

It is also useful to make the following remark. The moduli space of classical vacua $\Ms_0$ is well known to have several irreducible components, typically referred to as ``branches''. This feature was considered in \cite{FHHZ08a} to be reflected by the fact that $\uRep(A,\alpha)$ is in general a reducible variety. In this reference it was argued that there exists a unique top-dimensional irreducible component of $\uRep(A,\alpha)$, which we denote by $\Vs_A(\alpha)$. The ``Higgs branch'' of the moduli space of classical vacua is given by the affine quotient $\Vs_A(\alpha) \sslash \uGL(\alpha)$. The other components of $\Ms_0$ are commonly referred to as ``Coulomb branches''. 

Now we come to the central point. One of the main insights in \cite{KW99}, further explored in \cite{BGLP00}, was that for the $N=1$ quiver gauge theory associated to a D$3$-brane on an isolated Calabi-Yau singularity $X$, the Higgs branch of the moduli space of classical vacua recovers the geometry of $X$. A bit more precisely, we are asserting that
$$
X \cong \Vs_A(\alpha) \sslash \uGL(\alpha) \subseteq \Ms_0.
$$
Here the dimension vector $\alpha$ is fixed in terms of the multiplicities $\alpha_i$ of the fractional branes. From this perspective, the worldvolume gauge theory of the D-brane is the primary concept, whereas the spacetime itself is a secondary, derived concept.

The above discussion suggests that it would be possible to reconstruct singular Calabi-Yau threefolds from $N=1$ quiver gauge theories. Following the terminology in \cite{Ber02}, we call this process ``reverse geometric engineering of singularities''.  The basic idea of this construction may be summarised as follows. One is given an $N=1$ quiver gauge theory so that the corresponding $N=1$ quiver algebra $A$ is finitely generated as a module over its centre $Z(A)$. Then $Z(A)$ is itself the coordinate ring of a three-dimensional variety, to be identified as the singularity $X$. We must show, then, that the Higgs branch of the moduli space of classical vacua coincides with the variety $\uSpec Z(A)$. (This ties in with a general principle of noncommutative algebraic geometry espoused in \cite{Kont93} and further developed in \cite{Br06}.) A more physical version of this statement is to say that $\uSpec Z(A)$ will correspond to the ``spacetime'' in which closed strings propagate, while $A$ is associated to a noncommutative algebraic geometry that D-branes see. 

In a similar vein, one can expect that crepant resolutions of the Calabi-Yau singularity $X=\uSpec Z(A)$ are realized as moduli spaces of the D-brane theory, in the presence of Fayet-Iliopoulos parameters. Let us spell out more clearly what we mean by this. If $\theta\colon \Zx^I \fl \Zx$ satisfies $\theta(\alpha)=0$, then there are notions of $\theta$-stable and $\theta$-semistable elements of $\Vs_A(\alpha)$, there is a GIT quotient $\Vs_A(\alpha) \sslash_{\theta}\uGL(\alpha)$, and a natural map
$$
\pi_{\theta}\colon \Vs_A(\alpha) \sslash_{\theta}\uGL(\alpha) \fl \Vs_A(\alpha) \sslash \uGL(\alpha)
$$
which is a projective morphism; see Sect.~\ref{sec2.3} below for a more precise description. Moreover, if $\alpha$ is such that the general element of $\Vs_A(\alpha)$ is a simple representation of $A$, then $\pi_{\theta}$ is a birational map of irreducible varieties. The quotient $\Vs_A(\alpha) \sslash_{\theta}\uGL(\alpha)$ is a quasiprojective variety whose points are in bijection with S-equivalence classes of $\theta$-semistable elements of $\Vs_A(\alpha)$. This is usually called the ``mesonic moduli space with Fayet-Iliopoulos parameters $\theta$''. It turns out that, in many cases, if $\theta$ is chosen so that $\theta$-semistables are $\theta$-stable, then the corresponding mesonic moduli space $\Vs_A(\alpha) \sslash_{\theta}\uGL(\alpha)$ is a smooth Calabi-Yau threefold and so $\pi_{\theta}$ is a crepant resolution of the Calabi-Yau singularity $X=\uSpec Z(A)$.

We actually can go somewhat further along these lines. Let $Y$ be {\it any} crepant resolution of the singularities of $X=\uSpec Z(A)$. A point that we mentioned earlier but did not elaborate upon, is that D-branes on $Y$ should be properly regarded as objects in $\Dd^b(\uCoh(Y))$, the bounded derived category of coherent sheaves on $Y$. The foregoing discussion makes it highly plausible that the $N=1$ quiver gauge theory would give a different description of these D-branes in terms of representations of the $N=1$ quiver algebra $A$. This statement can be made more precise by saying that there is an equivalence of triangulated categories
$$
\Dd^b(\uCoh(Y)) \cong \Dd^b(\umod A),
$$
where $\Dd^b(\umod A)$ is the bounded derived category of finitely generated right modules over $A$. In M.~Van den Bergh's terminology, $A$ is a ``noncommutative crepant resolution'' of $X=\uSpec Z(A)$. 

The situation that we actually wish to apply this to is the case of ADE fibered Calabi-Yau threefolds and their small resolutions; see Sect.~\ref{setup} below for details. The relevant $N=1$ quiver gauge theory was written down in \cite{CKV01} (see also \cite{CFIKV02,K04,Zhu06}). We will explicitly carry out the previous construction in the subsequent sections.

\section{ADE fibered Calabi-Yau threefolds and their small resolutions revisited}\label{sec2}
This section studies the reverse geometric engineering of ADE fibered Calabi-Yau threefolds and their small resolutions along the lines indicated in the previous section. Before we do so, we describe our general setup and fix notation.

\subsection{General setup}\label{setup}
Let $G$ be a finite subgroup of $\uSL(2,\Cx)$, let $\Cx^2/G$ be the corresponding Kleinian singularity and let $\pi_0 \colon Y_0 \fl \Cx^2/G$ be its minimal resolution. The exceptional divisor $C$ of $\pi_0$ is known to be a union of projective lines intersecting transversally, and the graph $\Gamma$ whose vertices correspond to the irreducible components of $C$, with two vertices joined if and only if the components intersect, is a Dynkin diagram of type $A$, $D$ or $E$. 

Now let $Z(\Cx G)$ be the centre of the group algebra $\Cx G$, and let $\hfrak$ be the codimension one hyperplane in $Z(\Cx G)$ formed by all central elements which have trace zero in $\Cx G$. According to the McKay correspondence, the dual space $\hfrak^*$ carries a root system associated to the Dynkin diagram $\Gamma$. Write $\mathfrak{W}$ for the Weyl group of this root system. It is then a fairly standard result (see, for example, \cite{Slodow80} and references therein) that the Kleinian singularity $\Cx^2/G$ has a semi-universal deformation, a flat family $\varphi \colon\Xs \fl \hfrak/ \mathfrak{W}$ whose fiber over $0$ is $\Cx^2/G$.  Furthermore, Brieskorn and Tyurina \cite{Briesk68,Tj70} showed that this family admits a simultaneous resolution after making the base change $\mathfrak{h} \fl \mathfrak{h}/\mathfrak{W}$. More precisely, the family $\Xs \times_{\mathfrak{h}/\mathfrak{W}} \mathfrak{h}$ may be resolved explicitly and one obtains a simultaneous resolution $\Ys \fl \Xs \times_{\mathfrak{h}/\mathfrak{W}} \mathfrak{h}$ of $\varphi$ inducing the minimal resolution $Y_0 \fl \Cx^2/G$. The situation can be conveniently summarised by the diagram
$$
\SelectTips{cm}{10}\xymatrix{\Ys \ar[r] \ar[rd] &\Xs \times_{\mathfrak{h}/\mathfrak{W}} \mathfrak{h} \ar[r] \ar[d] & \Xs \ar[d] \\
                                   & \mathfrak{h} \ar[r] &  \mathfrak{h}/\mathfrak{W}}
$$

Using these observations we can define a broader class of Calabi-Yau threefolds as follows. We want to obtain a Gorenstein Calabi-Yau threefold $X$ by fibering the total space of the semiuniversal deformation of a Kleinian singularity $\Cx^2/G$ over a complex plane. To make things more concrete, let $t\colon  \Cx \fl \hfrak$ be a polynomial map. Via the defining equation for the family $\Xs \times_{\mathfrak{h}/\mathfrak{W}} \mathfrak{h}$, we can view $X$ as the total space of a one parameter family defined by $t$. Similarly, the simultaneous resolution $\Ys \fl \Xs \times_{\mathfrak{h}/\mathfrak{W}} \mathfrak{h}$ can be used to construct a Calabi-Yau threefold $Y$. That is, we get a cartesian diagram
$$
\SelectTips{cm}{10}\xymatrix{Y \ar[r] \ar[d]_-{\pi} & \Ys \ar[d] \\
              X \ar[r] \ar[d] & \Xs \times_{\mathfrak{h}/\mathfrak{W}} \mathfrak{h} \ar[d] \\
              \Cx \ar[r]^-{t} & \mathfrak{h}}
$$
where $Y$ is the pullback of $\Ys$ by $t$ and $X$ is the pullback of $\Xs \times_{\mathfrak{h}/\mathfrak{W}} \mathfrak{h}$ by $t$. One can show that if $t$ is sufficiently general, then $Y$ is smooth, $X$ is Gorenstein with an isolated singular point, and $\pi \colon  Y \fl X$ is a small resolution. The genericity condition is that $t$ is transverse to the hyperplanes $\Pi_{\rho} \subseteq \hfrak$ orthogonal to each positive root $\rho$ of $\hfrak$. This class of Calabi-Yau threefolds we call ADE fibered Calabi-Yau threefolds. 

Note that if $t=0$, then the singular threefold $X$ is isomorphic to a direct product of the form $\Cx^2/G \times \Cx$. In particular $X$ has a line of Kleinian singularities. The resolution $Y$ is isomorphic to the direct product $Y_0 \times \Cx$. This of course is {\it not} a small resolution, as it has an exceptional divisor over a curve. The main point for us here is that ADE fibered Calabi-Yau threefolds are related to $\Cx^2/G \times \Cx$ by a complex structure deformation. For a full discussion of these matters consult \cite{Sz04,Sz08}.   

\subsection{ADE fibered Calabi-Yau threefolds revisited}\label{ADEquiverAlg}
In this subsection we show how to construct ADE fibered Calabi-Yau threefolds in terms of a noncommutative algebra, which we call the ``$N=1$ ADE quiver algebra''. This confirms what was suggested earlier more informally. We start by summarising some of the neccessary definitions.

Let $Q$ be an extended Dynkin quiver with vertex set $I$, and let $h(a)$ and $t(a)$ denote the head and tail vertices of an arrow $a \in Q$. The double $\overline{Q}$ of $Q$ is the quiver obtained by adding a reverse arrow $a^*\colon j \fl i$ for each arrow $a \colon i \fl j$ in $Q$. We denote by $\widehat{Q}$ the quiver obtained from $\overline{Q}$ by attaching an additional edge-loop $u_i$ for each vertex $i \in I$. We write $\Cx \overline{Q}$ and $\Cx \widehat{Q}$ for the path algebras of $\overline{Q}$ and $\widehat{Q}$. 

Now let $B=\bigoplus_{i \in I} \Cx e_i$ be the semisimple commutative subalgebra of $\Cx \overline{Q}$ spanned by the trivial paths and consider the algebra $B[u]$ of polynomials in an indeterminate $u$ with coefficients in $B$. For an element $\tau \in B[u]$, we will write $\tau(u)=\sum_{i}\tau_i(u) e_i$ where $\tau_i \in \Cx[u]$. Let $\Cx \overline{Q} \ast_B B[u]$ denote the free product\footnote{``Free product'' is a misnomer in this context: if $R$ is a commutative ring and $A$ and $B$ are $R$-algebras then $A \ast_R B$ is the coproduct in the category of $R$-algebras.} of $\Cx \overline{Q}$ with $B[u]$ over $B$. We have an isomorphism
$$
\SelectTips{cm}{10}\xymatrix{\Cx \overline{Q} \ast_B B[u]~\widetilde{\lfl}~\Cx \widehat{Q} \colon  u \mapstochar\longrightarrow \sum_i u_i.}
$$
This isomorphism sends the element $e_i u e_i$ to $u_i$, the additional edge-loop at the vertex $i$. We also have an isomorphism
$$
\SelectTips{cm}{10}\xymatrix{B[u]= \Big(\bigoplus_{i \in I}\Cx e_i \Big) \otimes\Cx[u]~\widetilde{\lfl}~\bigoplus_{i \in I}\Cx[u_i]\colon  e_i \otimes u \mapstochar\longrightarrow u_i.}
$$
Therefore, choosing an element $\tau \in B[u]$ amounts to choosing a collection of polynomials $\{ \tau_i \in \Cx[u_i] \mid i \in I \}$. 

If $\tau \in B[u]$ then the $N=1$ {\it ADE quiver algebra} determined by $Q$ is defined by
$$
\mathfrak{A}^{\tau}(Q)=\Cx \widehat{Q} \bigg/ \left( \begin{array}{c} \displaystyle\sum\nolimits_{a \in Q}[a,a^*]-\displaystyle\sum\nolimits_{i \in I} \tau_i(u) e_i \\ u \mbox{ {\small is a central element}} \end{array} \right).
$$
Compare this with the definition of \cite[\S 4.3]{Gin06}. We point out that the defining relations for $\mathfrak{A}^{\tau}(Q)$ are generated by the superpotential
$$
W=u \sum_{a \in Q}[a,a^*]-\sum_{i \in I}\eta_i(u) e_i,
$$
where each $\eta_i \in \Cx[u]$ satisfies $\eta'_i(u)=\tau_i(u)$. We also note that if $\tau(u)$ is identified with the element $\sum_i \tau_i(u_i)$ then $\mathfrak{A}^{\tau}(Q)$ is the same as the quotient of $\Cx \widehat{Q}$ by the relations
$$
\sum_{h(a)=i}aa^*-\sum_{t(a)=i}a^*a-\tau_i(u_i)=0, \quad a u_i=u_j a,
$$
for each vertex $i$, and for each arrow $a\colon i \fl j$ in $\overline{Q}$. This is helpful when considering representations of $\mathfrak{A}^{\tau}(Q)$, as they can be identified with representations $V$ of $\widehat{Q}$ which satisfy
\begin{equation*}\label{eqnREL1}
\tag{\textasteriskcentered}
\sum_{h(a)=i}V_aV_{a^*}-\sum_{ t(a)=i}V_{a^*}V_a-\tau_i(V_{u_i})=0, \quad V_a V_{u_i}=V_{u_j} V_a,
\end{equation*}
for each vertex $i$, and for each arrow $a\colon i \fl j$ in $\overline{Q}$.   

The following is immediate from what we have just seen.

\begin{lem}\label{lem2.1}
Let $V$ be a simple representation of $\mathfrak{A}^{\tau}(Q)$. Then there exists a $\lambda$ such that $V_{u_i} v=\lambda v$ for all $i \in I$ and all $v \in V_i$. 
\end{lem}  

\begin{proof}
By virtue of (\ref{eqnREL1}), we may regard the collection of $\Cx$-linear maps $\{V_{u_i}\colon V_i \fl V_i \mid i \in I \}$ as an endomorphism $\phi\colon V \fl V$. Because $V$ is simple, it follows from Schur's lemma that $\phi$ must act as a scalar multiple $\lambda \cdot \id_V$ of the identity. The latter obviously implies the assertion.
\end{proof}

If $\alpha \in \Nx^I$, then representations of $\widehat{Q}$ of dimension vector $\alpha$ are given by elements of the variety
$$
\uRep(\widehat{Q},\alpha)=\Bigg(\bigoplus_{a\in Q}\uHom_{\Cx}(\Cx^{\alpha_{t(a)}},\Cx^{\alpha_{h(a)}}) \oplus \uHom_{\Cx}(\Cx^{\alpha_{h(a)}},\Cx^{\alpha_{t(a)}}) \Bigg) \oplus \Bigg(\bigoplus_{i \in I} \uEnd_{\Cx}(\Cx^{\alpha_i})\Bigg).
$$
We denote by $\uRep(\mathfrak{A}^{\tau}(Q),\alpha)$ the closed subspace of $\uRep(\widehat{Q},\alpha)$ corresponding to representations for $\mathfrak{A}^{\tau}(Q)$. The group $\uGL(\alpha)$ acts on both these spaces, and the orbits correspond to isomorphism classes. 
 
 We have the following easily verified result.

\begin{lem}\label{lem2.2}
If $x \in  \uRep(\mathfrak{A}^{\tau}(Q),\alpha)$ then $\sum_i \utr \tau_i(x_{u_i})=0$.
\end{lem}

\begin{proof}
Given $a \in \overline{Q}$, we have $\utr(x_a x_{a^*})=\utr(x_{a^*}x_{a})$. Taking traces to relations (\ref{eqnREL1}) and summing over all vertices $i \in I$, one obtains $\sum_i \utr \tau_i(x_{u_i})=0$, as required.
\end{proof}

We define $\Vs_Q(\tau,\alpha)$ to be the subset of $\uRep(\mathfrak{A}^{\tau}(Q),\alpha)$ consisting of the representations $x$ for which there exists a $\lambda$ (depending on $x$) with the property that each $x_{u_i}$ acts diagonally on $\Cx^{\alpha_i}$ by multiplication by $\lambda$, i.e., for any $i \in I$ and $v \in V_i$, we will have $x_{u_i}v=\lambda v$. It is clear that this is a locally closed subset of $\uRep(\mathfrak{A}^{\tau}(Q),\alpha)$, so a variety. In view of Lemma~\ref{lem2.1}, we immediately deduce that $\Vs_Q(\tau,\alpha)$ contains the open subset $\Ss_Q(\tau,\alpha)$ consisting of simple representations of $\mathfrak{A}^{\tau}(Q)$. 

The next result is an easy consequence of Lemma~\ref{lem2.2}.

\begin{cor}\label{cor2.3}
If $x \in \Vs_Q(\tau,\alpha)$ then $\sum_i \alpha_i \tau_i(\lambda)=0$.
\end{cor} 

We are now ready to start our study of ADE fibered Calabi-Yau threefolds. Let $G$ be a finite group of $\uSL(2,\Cx)$ and let $\Cx^2/G$ be the corresponding Kleinian singularity. Let $\rho_0,\dots,\rho_n$ be the irreducible representations of $G$ with $\rho_0$ trivial, and let $V$ be the natural $2$-dimensional representation of $G$. The McKay graph of $G$ is the graph with vertex set $I=\{0,1,\dots,n \}$ and with the number of edges between $i$ and $j$ being the multiplicity of $\rho_i$ in $V \otimes \rho_j$. According to the McKay correspondence, this graph is an extended Dynkin diagram of type $A$, $D$ or $E$. Let $Q$ be the quiver obtained from the McKay graph by choosing any orientation of the edges, and let $\delta \in \Nx^I$ be the vector with $\delta_i = \dim \rho_i$. 

In \cite{CKV01}, Sect.~5, it is pointed out that the moduli space of classical vacua $\Ms_0 = \uRep(\mathfrak{A}^{\tau}(Q),\delta) \sslash \uGL(\delta)$ contains a branch (the Higgs branch) which is an ADE fibered Calabi-Yau threefold. With the above formulation, we wish to present a clear proof of this claim. First, however, it will be convenient to provide the following piece of information. The space of representations of $\overline{Q}$ of dimension vector $\delta$ can be identified with a cotangent bundle 
$$
\uRep(\overline{Q},\delta)\cong \uRep(Q,\delta) \times \uRep(Q,\delta)^* \cong T^* \uRep(Q,\delta).
$$
This has a natural symplectic structure, and associated to the action of $\uGL(\delta)$ there is a moment map $\mu_{\delta}\colon \uRep(\overline{Q},\delta) \fl \bigoplus_{i \in I}\mathfrak{gl}(\delta_i,\Cx)$ given by
$$
\mu_{\delta}(x)_i=\sum_{h(a)=i}x_a x_{a^*}-\sum_{t(a)=i}x_{a^*}x_{a}.
$$
Identifying $\tau \in \Cx^I$ with the element of $\bigoplus_{i \in I}\mathfrak{gl}(\delta_i,\Cx)$ whose $i$th component is $\tau_i$ times the identity matrix of size $\delta_i$, one can consider the fiber $\mu_{\delta}^{-1}(\tau)$ and the affine quotient variety $\mu_{\delta}^{-1}(\tau) \sslash \uGL(\delta)$. The elements of $\mu_{\delta}^{-1}(\tau)$ correspond to representations of dimension vector $\delta$ of a certain algebra, the {\it deformed preprojective algebra} $\Pi^{\tau}(Q)$ of \cite{CrH98}. More directly relevant for us is that, letting $\hfrak  = \{ \tau \in \Cx^I \mid \delta \cdot \tau =0\}$, there is a flat family $\varphi\colon \mu_{\delta}^{-1}(\hfrak) \sslash \uGL(\delta) \fl \hfrak$ whose fiber over $0$ is $\Cx^2/G$. By the discussion at the end of Sect.~8 of \cite{CrH98}, this family is obtained from the semi-universal deformation of $\Cx^2/G$ by lifting through the Weyl group. 

Now let $x$ be any element of $\Vs_Q(\tau,\delta)$. By Corollary~\ref{cor2.3}, the dimension vector $\delta$ satisfies $\sum_i \delta_i \tau_i(\lambda)=0$. Furthermore, for any vertex $i$ we have $\mu_{\delta}(x)_i =\tau_i(\lambda)$. Therefore one can identify $\Vs_Q(\tau,\delta)$ with the fiber product
$$
\SelectTips{cm}{10}\xymatrix{\Vs_Q(\tau,\delta) \ar[r] \ar[d] & \mu_{\delta}^{-1}(\hfrak) \ar[d]^-{\mu_{\delta}}\\
              \Cx \ar[r]^-{\widehat{\tau}} & \hfrak}
$$
where $\widehat{\tau}\colon  \Cx \fl \hfrak$ is the map corresponding to $\tau$. Observe that $\uGL(\delta)$ acts naturally on $\Vs_Q(\tau,\delta)$ in such a way that all maps in the fiber product are equivariant (where the action on $\Cx$ is trivial). Now $\Vs_Q(\tau,\delta) \fl \Cx$ is flat since it is the pullback of $\mu_{\delta}$, which is flat by \cite[Lemma~8.3]{CrH98}. From this it follows that the map $\Vs_Q(\tau,\delta)\sslash \uGL(\delta) \fl \Cx$ is also flat and surjective.

Incidentally, if $\tau \in \hfrak$, then $\mu_{\delta}^{-1}(\tau)$ is irreducible by \cite[Lemma~6.3]{Cr01}, which implies that every fiber of the map $\Vs_Q(\tau,\delta) \fl \Cx$ is irreducible. It then follows from \cite[Lemma~6.1]{Cr01} that $\Vs_Q(\tau,\delta)$ is irreducible. On the other hand, we may apply Theorem~6.7 of \cite{Cr01} to infer that $\mathfrak{A}^{\tau}(Q)$ has a simple representation of dimension vector $\delta$. This allows us to conclude that the set $\Ss_Q(\tau,\delta)$ of simple representations must be dense.
 
We have now accumulated all the information necessary to prove the following result.

\begin{thm}\label{thmNC3fold}
Assume that $\widehat{\tau}$ is sufficiently general. Then the affine quotient variety $\Vs_Q(\tau,\delta)\sslash \uGL(\delta)$ is isomorphic to the ADE fibered Calabi-Yau threefold associated with $Q$ and $\widehat{\tau}$.
\end{thm}

\begin{proof}
Let $\Cx[\mu_{\delta}^{-1}(\hfrak)]$ and $\Cx[\hfrak]$ be the coordinate rings of $\mu_{\delta}^{-1}(\hfrak)$ and $\hfrak$, respectively. Then the coordinate ring of $\Vs_Q(\tau,\delta)$ is given by
$$
\Cx[\Vs_Q(\tau,\delta)] = \Cx[\mu_{\delta}^{-1}(\hfrak)] \otimes_{\Cx[\hfrak]}\Cx[u].
$$
Since $\uGL(\delta)$ is linear reductive and acts trivially on $\Cx[u]$, we see that
$$
\Cx[\Vs_Q(\tau,\delta)]^{\uGL(\delta)}  \cong \Cx[\mu_{\delta}^{-1}(\hfrak)]^{\uGL(\delta)} \otimes_{\Cx[\hfrak]} \Cx[u].
$$
Accordingly, we have $\Vs_Q(\tau,\delta)\sslash \uGL(\delta) \cong \mu_{\delta}^{-1}(\hfrak) \sslash \uGL(\delta) \times_{\hfrak} \Cx$. Hence we obtain the affine quotient $\Vs_Q(\tau,\delta)\sslash \uGL(\delta)$ as the fiber product
$$
\SelectTips{cm}{10}\xymatrix{\Vs_Q(\tau,\delta)\sslash \uGL(\delta) \ar[r] \ar[d] & \mu_{\delta}^{-1}(\hfrak)\sslash \uGL(\delta) \ar[d]^-{\varphi}\\
              \Cx \ar[r]^-{\widehat{\tau}} & \hfrak}
$$
Since, by hypothesis, $\widehat{\tau}$ is sufficiently general, the desired assertion follows. 
\end{proof} 

We illustrate with the following concrete example.

\begin{exam}\label{exam3.6}
Suppose that $Q$ is of type $\widetilde{A}_n$, so that $\delta_i=1$ for all vertices $i$. The arrows $a_i$ in $Q$ connect vertices $i$ and $i+1$ (identifying $n+1$ with zero). Thus $\widehat{Q}$ has shape
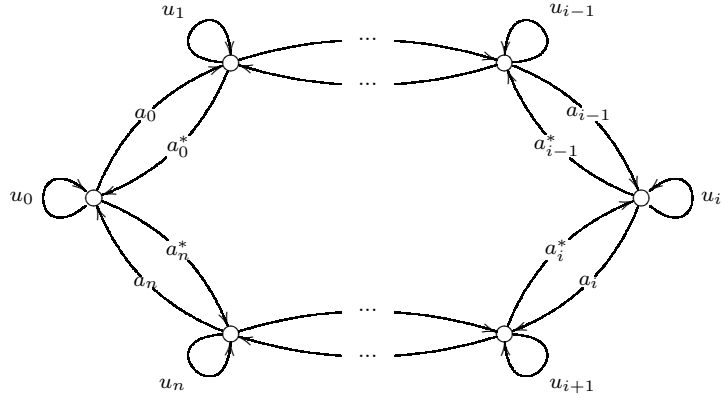
\begin{equation*}
\begin{xy} <1.0mm,0mm>:
(0,0)*{\cir<3pt>{}}="a",(18,18)*{\cir<3pt>{}}="b",(54,18)*{\cir<3pt>{}}="c",(72,0)*{\cir<3pt>{}}="d",(54,-18)*{\cir<3pt>{}}="e",(18,-18)*{\cir<3pt>{}}="f"
\ar@{->}@/^3mm/|{a_0}"a";"b"
\ar@{->}@/^3mm/|{\hole \dots  \hole}"b";"c"
\ar@{->}@/^3mm/|{a_{i-1}}"c";"d"
\ar@{->}@/^3mm/|{a_i}"d";"e"
\ar@{->}@/^3mm/|{\hole \dots  \hole}"e";"f"
\ar@{->}@/^3mm/|{a_n}"f";"a"
\ar@{<-}@/_3mm/|{a_0^*}"a";"b"
\ar@{<-}@/_3mm/|{\hole \dots  \hole}"b";"c"
\ar@{<-}@/_3mm/|{a_{i-1}^*}"c";"d"
\ar@{<-}@/_3mm/|{a_i^*}"d";"e"
\ar@{<-}@/_3mm/|{\hole \dots  \hole}"e";"f"
\ar@{<-}@/_3mm/|{a_n^*}"f";"a"
\ar@(dl,ul)^-{u_0}"a";"a"
\ar@(l,u)^-{u_1}"b";"b"
\ar@(r,u)_-{u_{i-1}}"c";"c"
\ar@(dr,ur)_-{u_{i}}"d";"d"
\ar@(r,d)^-{u_{i+1}}"e";"e"
\ar@(l,d)_-{u_n}"f";"f"
\end{xy}
\end{equation*}
As above let $\tau \in B[u]$; recall that it is specified by a set of polynomials $\{ \tau_i \in \Cx[u_i] \mid 0 \leq i \leq n \}$. Because $\delta=(1,\dots,1)$, a representation of $\widehat{Q}$ of dimension $\delta$ involves placing a one-dimensional vector space at each vertex $i$ and assigning a complex number to each arrow $a_i,a^*_i,u_i$. Hence, we may identify $\uRep(\widehat{Q},\delta)$ with the space $\Cx^{n+1} \times \Cx^{n+1} \times \Cx^{n+1}$ so that
\begin{align*}
\uRep(\mathfrak{A}^{\tau}(Q),\delta) & \cong \{ (x_i,y_i,\lambda_i) \mid -x_i y_i + x_{i+1}y_{i+1} = \tau_i(\lambda_i), 0 \leq i \leq n \} \\
& \subseteq \Cx^{n+1} \times \Cx^{n+1} \times \Cx^{n+1}.
\end{align*}
Also one can identify
\begin{align*}
\Vs_Q(\tau,\delta) & \cong \{ (x_i,y_i,\lambda) \mid  -x_i y_i + x_{i+1}y_{i+1} = \tau_i(\lambda), 0 \leq i \leq n  \} \\
&\subseteq \Cx^{n+1} \times \Cx^{n+1} \times \Cx.
\end{align*}
The relations for $\Vs_Q(\tau,\delta)$ lead to the condition that $\sum_{i=0}^n \tau_i(\lambda)=0$. From this it follows that the map $\widehat{\tau}=(\tau_0,\dots,\tau_n)\colon \Cx \fl \Cx^{n+1}$ corresponding to $\tau$ has its image in $\hfrak$. Without loss of generality it is possible to suppose that $\tau_i=t_i-t_{i+1}$ for some polynomial map $t=(t_0,\dots,t_n)\colon \Cx \fl \hfrak$.

Now the action of $\uGL(\delta)=(\Cx^{\times})^{n+1}$ on $\Vs_Q(\tau,\delta)$ is by
$$
(x_i,y_i,\lambda)  \mapstochar\longrightarrow (g_{i+1} g_{i}^{-1}x_i,g_{i} g_{i+1}^{-1}y_i,\lambda)
$$ 
for $(g_i) \in \uGL(\delta)$ and $(x_i,y_i,\lambda) \in \Vs_Q(\tau,\delta)$. It is easily seen that the ring of invariants $\Cx[\Vs_Q(\tau,\delta)]^{\uGL(\delta)}$ is generated by
\begin{align*}
x & = x_0 \cdots x_n,\\
y &= y_0 \cdots y_n,\\
z_i &= x_i y_i, \quad 0 \leq i \leq n.
\end{align*}
These invariants satisfy the relation
$$
xy=z_0 \cdots z_n.
$$
On the other hand, the relations for $\Vs_Q(\tau,\delta)$ imply that
$$
z_i= z_n - \sum_{j=0}^i \tau_j(\lambda), \quad 0 \leq i \leq n.
$$
Bearing in mind that $\sum_{j=0}^i \tau_j(\lambda)=t_0(\lambda)-t_{i+1}(\lambda)$, we derive 
$$
\sum_{i=0}^n z_i =(n+1) (z_n-t_0(\lambda)).  
$$
Setting $z=\frac{1}{n+1}\sum_{i=0}^n z_i$, we therefore deduce that
$$
z_i = z+ t_{i+1}(\lambda), \quad 0 \leq i \leq n.
$$
The conclusion is that the affine quotient variety $\Vs_Q(\tau,\delta)\sslash \uGL(\delta)$ is given as the hypersurface
$$
\bigg\{ (x,y,z,\lambda) \,\bigg|\, xy=\prod_{i=0}^n (z + t_{i+1}(\lambda)) \bigg\} \subseteq \Cx^4
$$
which is the total space of the family describing the $A_n$ fibration over the $\lambda$-plane.
\end{exam}

\subsection{Small resolutions of ADE fibered Calabi-Yau threefolds revisited}\label{sec2.3}
We continue with the same hypothesis and notation as in the previous subsection. Our main aim in this paragraph is to show how small resolutions of ADE fibered Calabi-Yau threefolds can be obtained as moduli spaces of representations of the $N=1$ ADE quiver algebra $\mathfrak{A}^{\tau}(Q)$. This is a further confirmation of the picture we have developed. We begin with some generalities.

Let $A$ be a commutative ring, let $S$ be a graded ring with $S_0=A$, and let $\phi\colon A \fl R$ be a commutative ring homomorphism. Then we may construct the tensor product $T=S \otimes_A R$ of $S$ and $R$ over $A$ by means of $\phi$. We consider $T$ as a graded ring by $T_d=S_d \otimes_A R$. Clearly $T_0=R$. We then get the following simple observation.

\begin{lem}\label{lem2.6}
With the notation above, we have $\uProj T \cong \uProj S \times_{\uSpec A}\uSpec R$.
\end{lem}

\begin{proof}
We note first that $\phi \colon A \fl R$ extends to a homomorphism $\psi \colon S \fl T$ of graded rings (preserving degrees). Let $U = \{ \mathfrak{p} \in \uProj T \mid \mathfrak{p} \nsupseteq \psi(S_+) \}$ and $\widetilde{\psi} \colon U \fl \uProj S$ be the morphism determined by $\psi$. Since $\psi(S_+)$ generates $T_+$, it is clear that $U = \uProj T$. Now fix a homogeneous $f \in S$, and set $g = \psi(f)$. Then we get an isomorphism from $S_{(f)} \otimes_A R$ to $T_{(g)}$ by assigning $s/f^n \otimes r$ to $(s \otimes r)/g^n$. This implies that $\uSpec T_{(g)} \cong \uSpec S_{(f)} \times_{\uSpec A}\uSpec R$, whence the result follows.
\end{proof}

Let us now proceed with the construction of small resolutions of ADE fibered Calabi-Yau threefolds. For convenience of reference, we first record some background information. Recall that the flat family $\varphi \colon \mu_{\delta}^{-1}(\hfrak) \sslash \uGL(\delta) \fl \hfrak$ realizes the semi-universal deformation of $\Cx^2 / G$, or rather its lift through a Weyl group action. Let $\theta \colon \Zx^I \fl \Zx$ satisfy $\theta(\delta)=0$. Then $\theta$ determines a character $\chi_{\theta}$ of $\uGL(\delta)$ mapping $(g_i)$ to $\prod_i \det (g_i)^{\theta_i}$. We denote by $\Cx[\mu_{\delta}^{-1}(\hfrak)]^{\uGL(\delta)}_{n\theta}$ the set of $\theta$-semi-invariants of weight $n$, i.e.~those functions on which $\uGL(\delta)$ acts by the character $\chi_{\theta}^n$. We say that $x \in \mu_{\delta}^{-1}(\hfrak)$ is $\theta$-semistable if there is a $\theta$-semi-invariant $f \in \Cx[\mu_{\delta}^{-1}(\hfrak)]^{\uGL(\delta)}_{n\theta}$ with $n \geq 1$, such that $f(x) \neq 0$. We say that $x$ is $\theta$-stable if, in addition, the stablizer of $x$ in $\uGL(\delta)$ is finite and the action of $\uGL(\delta)$ on $\{ x \in \mu_{\delta}^{-1}(\hfrak) \mid f(x)\neq 0 \}$ is closed. Following \cite{K94}, the GIT quotient $\mu_{\delta}^{-1}(\hfrak) \sslash_{\theta} \uGL(\delta)=\uProj \big(\bigoplus_{n \geq 0}\Cx[\mu_{\delta}^{-1}(\hfrak)]^{\uGL(\delta)}_{n\theta}\big)$ is the categorical quotient $\mu_{\delta}^{-1}(\hfrak)^{\mathrm{SS}}_{\theta}/ \uGL(\delta)$, where $\mu_{\delta}^{-1}(\hfrak)^{\mathrm{SS}}_{\theta}$ is the open subset parametrizing $\theta$-semistable elements of $\mu_{\delta}^{-1}(\hfrak)$. A parameter $\theta$ is generic if every $\theta$-semistable element is $\theta$-stable. In this case, $\mu_{\delta}^{-1}(\hfrak) \sslash_{\theta} \uGL(\delta)$ is the geometric quotient $\mu_{\delta}^{-1}(\hfrak)^{\mathrm{S}}_{\theta}\sslash \uGL(\delta)$, where $\mu_{\delta}^{-1}(\hfrak)^{\mathrm{S}}_{\theta}$ parametrizes $\theta$-stable elements of $\mu_{\delta}^{-1}(\hfrak)$. 
Since $\Cx[\mu_{\delta}^{-1}(\hfrak)]^{\uGL(\delta)}$ is a subalgebra of the graded ring defining $\mu_{\delta}^{-1}(\hfrak) \sslash_{\theta} \uGL(\delta)$, the $\uProj$ construction induces a projective morphism $\psi_{\theta}$ from $\mu_{\delta}^{-1}(\hfrak) \sslash_{\theta} \uGL(\delta)$ to $\mu_{\delta}^{-1}(\hfrak) \sslash \uGL(\delta)$. We write $\varphi_{\theta}$ to denote the composition of $\psi_{\theta}$ and $\varphi$. Using results of Cassens and Slodowy \cite[\S 7]{CS98} one can show that, for generic $\theta$, the family $\varphi_{\theta} \colon \mu_{\delta}^{-1}(\hfrak) \sslash_{\theta}  \uGL(\delta) \fl \hfrak$ induces a simultaneous resolution of $\varphi$.

Now consider $\Vs_Q(\tau,\delta) \subseteq \uRep(\mathfrak{A}^{\tau}(Q),\delta)$. Just as for $\mu_{\delta}^{-1}(\hfrak)$, one can consider $\theta$-semistable and $\theta$-stable points of $\Vs_Q(\tau,\delta)$. Again one can form a GIT quotient $\Vs_Q(\tau,\delta)\sslash_{\theta}  \uGL(\delta)$ and we get a projective morphism
$$
\pi_{\theta}\colon  \Vs_Q(\tau,\delta)\sslash_{\theta} \uGL(\delta) \lfl \Vs_Q(\tau,\delta)\sslash \uGL(\delta). 
$$
We have seen earlier that $\Vs_Q(\tau,\delta)$ is irreducible, and the general element is a simple representation of $\mathfrak{A}^{\tau}(Q)$, hence $\theta$-stable. Thus the morphism $\pi_{\theta}$ is a birational map of irreducible varieties.

We are at last in a position to attain our main objective, which is to prove the following result.

\begin{thm}\label{thm2.7}
Assume that $\widehat{\tau}$ is sufficiently general. If $\theta$ is generic, then $\pi_{\theta}$ is a small resolution of the ADE fibered Calabi-Yau threefold associated with $Q$ and $\widehat{\tau}$. 
\end{thm}

\begin{proof}
We keep the notation employed in the proof of Theorem~\ref{thmNC3fold}. To begin with, we observe that there is an isomorphism
$$
\Cx[\Vs_Q(\tau,\delta)]_{n\theta}^{\uGL(\delta)} \cong \Cx[\mu_{\delta}^{-1}(\hfrak)]_{n\theta}^{\uGL(\delta)} \otimes_{\Cx[\hfrak]} \Cx[u].
$$
Invoking Lemma~\ref{lem2.6}, it follows that
$$
\uProj  \Bigg(\bigoplus_{n \geq 0}\Cx[\Vs_Q(\tau,\delta)]_{n\theta}^{\uGL(\delta)} \Bigg) \cong \uProj  \Bigg(\bigoplus_{n \geq 0}\Cx[\mu_{\delta}^{-1}(\hfrak)]_{n\theta}^{\uGL(\delta)} \Bigg) \times_{\hfrak} \Cx,
$$
which entails $\Vs_Q(\tau,\delta)\sslash_{\theta} \uGL(\delta) \cong\mu_{\delta}^{-1}(\hfrak)\sslash_{\theta} \uGL(\delta) \times_{\hfrak} \Cx$. Therefore we obtain the GIT quotient $\Vs_Q(\tau,\delta)\sslash_{\theta} \uGL(\delta)$ as the fiber product
$$
\SelectTips{cm}{10}\xymatrix{\Vs_Q(\tau,\delta)\sslash_{\theta} \uGL(\delta) \ar[r] \ar[d] & \mu_{\delta}^{-1}(\hfrak)\sslash_{\theta} \uGL(\delta) \ar[d]^-{\varphi_{\theta}}\\
              \Cx \ar[r]^-{\widehat{\tau}} & \hfrak}
$$
The required result now follows from our hypothesis on $\widehat{\tau}$.
\end{proof}
 
We end this section with the following illustration of Theorem~\ref{thm2.7}.

\begin{exam}
Assume that $Q$ is a quiver of Dynkin type $\widetilde{A}_n$. We use the notations introduced in Example~\ref{exam3.6}. Consider the generic stability parameter $\theta=(-n,1,\dots,1)$. By definition, the ring of $\theta$-semi-invariants is spanned by the monomials $\prod_{i=0}^{n}x_i^{\alpha_i} y_{i}^{\beta_i}$ satisfying $-\alpha_0 + \alpha_n + \beta_0 - \beta_n =-n$ and $-\alpha_{i} + \alpha_{i-1} + \beta_i - \beta_{i-1}=1$ for $1 \leq i \leq n$. Given $j=0,\dots,n-1$, put
\begin{align*}
u_j&=x_0 \cdots x_j, \\
v_j&=y_{j+1} \cdots y_n.
\end{align*}
Then we have the following relations
\begin{align*}
& xv_j =u_j z_{j+1} \cdots z_n, \quad 0 \leq j \leq n-1, \\
& yu_j =v_j z_0 \cdots z_j, \quad 0 \leq j \leq n-1, \\
& u_j v_k = u_k v_j z_{k+1} \cdots z_j, \quad 0 \leq k < j \leq n-1,
\end{align*}
or, using the fact that $z_i=z+t_{i+1}(\lambda)$ for $0 \leq i \leq n$,
\begin{align*}
& xv_j=u_j \prod_{i=j+1}^n(z+t_{i+1}(\lambda)), \quad 0 \leq j \leq n-1, \\
& yu_j=v_j \prod_{i=0}^j(z+t_{i+1}(\lambda)), \quad 0 \leq j \leq n-1, \\
& u_j v_k=u_k v_j \prod_{i=k+1}^j(z+t_{i+1}(\lambda)), \quad 0 \leq k < j \leq n-1.
\end{align*}
Analyzing possibilities for $\alpha_i$, $\beta_i$ ($0 \leq i \leq n$) it is easily seen that the ring of $\theta$-semi-invariants is generated as a polynomial ring by 
$$
u_I v_{I'}=u_{i_1}\dots u_{i_p} v_{i'_1} \dots v_{i'_q},
$$
where $I=\{ i_1,\dots,i_p \}$ is a multi-index of $\{0,\dots, n-1 \}$ and $I'=\{ i'_1,\dots,i'_q \}$ denotes the complementary index. It is also not difficult to see that this space is the module over $\Cx[\Vs_Q(\tau,\delta)]^{\uGL(\delta)}$ generated by 
\begin{align*}
f_0 &=v_0 \cdots v_{n-2}v_{n-1}, \\
f_1 &=v_0 \cdots v_{n-2}u_{n-1}, \\
& \cdots \\
f_{n} &=u_0 \cdots u_{n-2}u_{n-1}.
\end{align*}
We observe next that
\begin{align*}
\Cx[u_I v_{I'} \mid I=\{ i_1,\dots,i_p \}, I'=\{ i'_1,\dots,i'_q \}] = \Cx[u_0,v_0] \ast \cdots \ast \Cx[u_{n-1},v_{n-1}],
\end{align*}
where $\ast$ denotes the Segre product of polynomial rings. Thus we have
\begin{align*}
\Cx[\Vs_Q&(\tau,\delta)]^{\uGL(\delta)}[f_0,\dots,f_n] \\
& \cong \Cx[\Vs_Q(\tau,\delta)]^{\uGL(\delta)}[u_0,v_0] \ast \cdots \ast \Cx[\Vs_Q(\tau,\delta)]^{\uGL(\delta)}[u_{n-1},v_{n-1}].
\end{align*}
Therefore the $\uProj$ quotient $\Vs_Q(\tau,\delta) \sslash_{\theta} \uGL(\delta)$ can be identified with a closed subvariety of $\Cx^{4} \times(\Px^1)^n$ with $(u_j : v_j)$ the homogeneous coordinates on the $j$th $\Px^1$. Now let $U_0, U_1,\dots, U_{n}$ be the open subsets of $\Cx^{4} \times(\Px^1)^n$ defined by
\begin{align*}
U_0 &=\{ v_0 \neq 0 \},\\
U_{k} &=\{ u_{k-1} \neq 0, v_{k} \neq 0 \}, \quad 1 \leq k \leq n-1,\\
U_{n} &=\{ u_{n-1} \neq 0 \},
\end{align*}
and on $U_k$, let
$$
\xi_k=v_{k-1}/u_{k-1}, \qquad \eta_{k}=u_k/v_k.
$$
Direct computations show that $\Vs_Q(\tau,\delta) \sslash_{\theta} \uGL(\delta) \cap U_k$ is defined by equations
\begin{align*}
(u_j:v_j)&= \bigg(1 : \xi_k \prod_{i=j+1}^{k-1}(z+t_{i+1}(\lambda)) \bigg), \quad \text{for $j<k-1$,} \\
(u_j:v_j)&=\bigg( \eta_{k} \prod_{i=k+1}^j(z+t_{i+1}(\lambda)) : 1\bigg), \quad \text{for $j>k$.}
\end{align*}
and
\begin{align*}
x &=\eta_{k}\prod_{i=k+1}^n (z + t_{i+1}(\lambda)), \\
 y &= \xi_k \prod_{i=0}^{k-1}(z + t_{i+1}(\lambda)), \\
z &=\xi_k \eta_k-t_{k+1}(\lambda).
\end{align*}
From the explicit analysis in \cite[\S 4]{KM92}, one sees that the GIT quotient $\Vs_Q(\tau,\delta) \sslash_{\theta} \uGL(\delta)$ is isomorphic to a small resolution of the threefold for a $A_n$ fibration over the $\lambda$-plane.   
\end{exam}

\section{Derived equivalence}\label{sec3}
In this section, it is shown how to describe the derived category of a small resolution of an ADE fibered Calabi-Yau threefold in terms of the associated $N=1$ ADE quiver algebra, in the spirit of noncommutative crepant resolutions of M.~Van den Bergh. Assertions of this sort have already been considered in \cite{Sz08}. Our work is mostly based on the ideas and constructions of \cite{CrH98} and \cite{Bergh02}.

\subsection{The algebra $A^{\tau}$}\label{subsec3.1}
As explained in the previous section, the $N=1$ ADE quiver algebra contains enough information to reconstruct ADE fibered Calabi-Yau threefolds and their small resolutions. Here we introduce yet another noncommutative algebra whose center is related to the coordinate ring of an ADE fibered Calabi-Yau threefold via Morita equivalence. The discussion borrows largely from \cite{CrH98}.

We begin by setting up notation. Let $G$ be a finite subgroup of $\uSL(2,\Cx)$. For our convenience, we denote by $R=\Cx[u]$ the ring of polynomials in a dummy variable $u$. The group $G$ acts naturally on the ring $R  \langle x,y\rangle$ of noncommuting polynomials, with the action of $G$ on $R$ being trivial, so one can form the skew group algebra $R \langle x,y\rangle \usm G$. We use $Z(RG)$ to denote the centre of the group algebra $RG$.

For $\tau \in Z(RG)$ we define the algebra $A^{\tau}$ as the quotient
$$
A^{\tau}=(R \langle x,y\rangle \usm G )/(xy-yx-\tau).
$$
This algebra was introduced and studied by W.~Crawley-Boevey and M.~Holland in \cite{CrH98}. (In the notation of \cite{CrH98}, $A^{\tau}$ corresponds to $\Ss^{R,\tau}$.) Observe that if $\tau=0$ then we recover the skew group algebra $R[x,y] \usm G$. In other words, $A^{\tau}$ is a flat deformation of $R[x,y] \usm G$ for every choice of $\tau$.  

The algebra $A^{\tau}$ carries a natural filtration, given by $\udeg x=\udeg y=1$, $\udeg u=0$ and $\udeg g=0$ for any $g \in G$. Let $\ugr A^{\tau}$ denote the associated graded algebra. It was explicitly demostrated in \cite[Lemma~1]{CrH98} that $\ugr A^{\tau} \cong R[x,y]\usm G$. As a consequence, the arguments in \cite[Sect.~1]{CrH98} go through and show that $A^{\tau}$ is a prime noetherian maximal order which is Auslander-Gorenstein and Cohen-Macaulay of GK dimension $3$.  The reader may want to consult \cite{StaZha94} for some background.

Now let $e=|G|^{-1}\sum_{g \in G}g$ be the averaging idempotent, viewed as an element in $A^{\tau}$. Define a subalgebra $C^{\tau}$ of $A^{\tau}$ to be $eA^{\tau}e$. The increasing filtration on $A^{\tau}$ induces a filtration on $C^{\tau}$. It is well known that $C^0 \cong R[x,y]^G$. Further $e$ lies in the degree zero part of the filtration of $A^{\tau}$ and therefore $\ugr C^{\tau} \cong e \ugr A^{\tau} e \cong R[x,y]^G$. This allows us to lift properties from $R[x,y]^G$ to $C^{\tau}$; we refer again to \cite{StaZha94} for details. In particular, $C^{\tau}$ is a noetherian integral domain of GK dimension $3$. Note finally that $C^{\tau}$ is a flat deformation of the coordinate ring of the associated singular Calabi-Yau threefold $\Cx^2/G \times \Cx$. 

We anticipate that if the trace of $\tau$ on the group algebra $RG$ is zero, then $C^{\tau}$ is a commutative ring, and it occurs as the coordinate ring of the fiber of the semi-universal deformation of $\Cx^2/G \times \Cx$. On the other hand, if $\tau$ has nonzero trace on $RG$, then $C^{\tau}$ is a noncommutative ring. 

In order to make further progress we need to bring in the notion of noncommutative crepant resolution introduced by M.~Van den Bergh \cite{Bergh02, Bergh04}. Let $R$ be an integrally closed Gorenstein domain. If $A$ is an $R$-algebra that is finite as an $R$-module, then $A$ is said to be homologically homogeneous if $A$ is a maximal Cohen-Macaulay $R$-module and $\ugldim A_{\mathfrak{p}}=\dim R_{\mathfrak{p}}$ for all $\mathfrak{p} \in \uSpec R$. A {\it noncommutative crepant resolution} of $R$ is a homologically homogeneous $R$-algebra of the form $A=\uEnd_{R}(M)$, where $M$ is a finitely generated reflexive $R$-module. We remind the reader that an $R$-module $M$ is said to be reflexive if the natural morphism $M \fl \uHom_{R}(\uHom_{R}(M,R),R)$ is an isomorphism.

Our aim now is to show that $A^{\tau}$ is a noncommutative crepant resolution of $C^{\tau}$. The following preliminary result will clear our path. 

\begin{lem}\label{lem3.1}
$A^{\tau} e$ is a finitely generated reflexive $C^{\tau}$-module. In addition, $A^{\tau} \cong \uEnd_{C^{\tau}}(A^{\tau}e)$.
\end{lem}

\begin{proof}
One can adapt the techniques of \cite[Lemma~1.4]{CrH98} to the present situation. Note first that $A^{\tau}eA^{\tau}$ is a finitely generated ideal of $A^{\tau}$. We write $A^{\tau}eA^{\tau}=\sum_{i=1}^m x_i A^{\tau}$ and $x_i=\sum_{j} r_{ij}e s_{ij}$ for some $r_{ij},s_{ij} \in A^{\tau}$. For each $a \in A^{\tau}$, we have $ae \in A^{\tau}e=(A^{\tau}eA^{\tau})e$, and so $ae = \big( \sum_i x_i a_i\big)e = \sum_{i,j}r_{ij} e s_{ij} a_i e$. This proves that the elements $r_{ij}$ generate $A^{\tau}e$ as a $C^{\tau}$-module. 

For the condition on the endomorphism ring, there are natural inclusions
$$
A^{\tau} \subseteq \uEnd_{C^{\tau}}(A^{\tau}e) \subseteq \uEnd_{A^{\tau}}(A^{\tau}eA^{\tau}).
$$
Let $Q$ denote the simple artinian quotient ring of $A^{\tau}$. The fact that $Q e \cong Q \otimes_{A^{\tau}}A^{\tau} e$ implies that $\uEnd_{C^{\tau}}(A^{\tau}e) \subseteq \uEnd_{C^{\tau}}(Qe)$. But $C^{\tau}$ is a maximal order in $eQe$, so $\uEnd_{C^{\tau}}(Qe)=\uEnd_{eQe}(Qe)$. Because $Q$ is simple artinian, we also have $Q \cong \uEnd_{eQe}(Qe)$. Thus the endomorphism ring $\uEnd_{C^{\tau}}(A^{\tau}e)$ can be identified with a subring of $Q$. From this it follows that
$$
\uEnd_{C^{\tau}}(A^{\tau}e) \cong \{ q \in Q \mid q A^{\tau}e \subseteq A^{\tau}e \}.
$$ 
Similarly, it can be shown that
$$
 \uEnd_{A^{\tau}}(A^{\tau}eA^{\tau}) \cong \{ q \in Q \mid q A^{\tau}eA^{\tau} \subseteq A^{\tau}eA^{\tau} \}=A^{\tau},
$$
the latter equality being an immediate consequence of the definition of a maximal order. The conclusion is that $A^{\tau} \cong \uEnd_{C^{\tau}}(A^{\tau}e)$, as asserted.

It remains to check that $A^{\tau}e$ is reflexive. A similar argument to the one above can be applied to show that $A^{\tau} \cong \uEnd_{C^{\tau}}(eA^{\tau})$. Hence 
$$
A^{\tau}e \cong \uHom_{C^{\tau}}(eA^{\tau}, eA^{\tau})e \cong \uHom_{C^{\tau}}(eA^{\tau}, C^{\tau}),
$$ 
and
$$
eA^{\tau} \cong e \uHom_{C^{\tau}}(A^{\tau}e, A^{\tau}e) \cong \uHom_{C^{\tau}}(A^{\tau}e, C^{\tau}),
$$
proving that $A^{\tau}e \cong \uHom_{C^{\tau}}(\uHom_{C^{\tau}}(A^{\tau}e, C^{\tau}), C^{\tau})$. This completes the proof of the lemma.
\end{proof}

We are now ready to prove our promised result.

\begin{prop}\label{prop3.2}
The algebra $A^{\tau}$ is a noncommutative crepant resolution of $C^{\tau}$.
\end{prop} 

\begin{proof}
By Lemma~\ref{lem3.1}, it suffices to show that $A^{\tau}$ is homologically homogeneous. We already know that $A^{\tau}$ is Cohen-Macaulay. Further, by Lemma~\ref{lem3.1} and \cite[Corollary~6.18]{McCR87}, $A^{\tau}$ has finite global dimension. The desired assertion now follows by appealing to \cite[Lemma~4.2]{Bergh02}.
\end{proof}

We now study the relationship between the algebra $A^{\tau}$ and the $N=1$ ADE quiver algebra. Keeping our earlier notation, the irreducible representations of $G$ are $\rho_0,\dots,\rho_n$, with $\rho_0$ trivial, and $I=\{0,1,\dots,n \}$. Let $Q$ be the quiver with vertex set $I$ obtained by choosing any orientation of the McKay graph, and let $\delta \in \Nx^I$ be the vector with $\delta_i=\dim \rho_i$. Fix an isomorphism $\Cx G \cong \bigoplus_{i \in I}\uMat_{\delta_i}(\Cx)$ and for every ordered pair $(p,q)$, $1 \leq p,q \leq \delta_i$, take $e_{ipq}$ to be the matrix with $p,q$ entry $1$ and zero elsewhere. Given $i \in I$, put $f_i = e_{i 1 1}$. Then $\{f_0,\dots, f_n \}$ is a set of nonzero orthogonal idempotents with the property $\Cx G f_i \cong \rho_i$ for all $i \in I$. Hence we get that $f=f_0 + \cdots +f_n$ is idempotent. Furthermore $f_0=e$, so $e=ef=fe$. Observe also that the map $R^I \fl Z(RG)$ given by $\tau \mapsto \sum_{i \in I}(\tau_i / \delta_i)f_i$ is a bijection, and we use this to identify $R^I$ and $Z(RG)$.

Before going on to give the connection between the algebra $A^{\tau}$ and the $N=1$ ADE quiver algebra $\mathfrak{A}^{\tau}(Q)$, it is convenient to point out the following description of $\mathfrak{A}^{\tau}(Q)$. We keep the notation of Sect.~\ref{ADEquiverAlg}. Following Crawley-Boevey and Holland \cite{CrH98},  given an element $\tau \in R^I$ we define $\Pi^{R,\tau}(Q)$ to be 
$$
R \overline{Q} \bigg/ \Bigg( \sum_{a \in Q}[a,a^*]-\sum_{i \in I} \tau_i e_i \Bigg).
$$
Because $u$ is central, we must have $R \overline{Q} = \Cx \overline{Q} \otimes R \cong \Cx \overline{Q} \ast_B B[u]$. Hence it follows that $\mathfrak{A}^{\tau}(Q) \cong \Pi^{R,\tau}(Q)$. With this understood, we get the following.

\begin{prop}\label{prop3.3}
$A^{\tau}$ is Morita equivalent to $\mathfrak{A}^{\tau}(Q)$ and $C^{\tau} \cong e_0 \mathfrak{A}^{\tau}(Q) e_0$.
\end{prop}

\begin{proof}
The first part of the proposition follows from the the fact that $f A^{\tau} f \cong \mathfrak{A}^{\tau}(Q)$ established in \cite[Theorem~3.4]{CrH98}. Under this isomorphism, $e$ corresponds to the trivial path $e_0$. Using that $e=ef=fe$, we have
$$
C^{\tau}=e A^{\tau} e = e f A^{\tau} f e \cong e_0 \mathfrak{A}^{\tau}(Q) e_0,
$$
as desired.
\end{proof}

For simplicity of notation we fix an isomorphism between $C^{\tau}$ and $e_0 \mathfrak{A}^{\tau}(Q) e_0$ and henceforth identify $C^{\tau} = e_0 \mathfrak{A}^{\tau}(Q) e_0$. Recall from Sect.~\ref{ADEquiverAlg} that one can identify $\Vs_Q(\tau,\delta)$ with the fiber product $ \mu_{\delta}^{-1}(\hfrak) \times_{\hfrak} \uSpec R$. Now, the coordinate ring of $\uRep(\overline{Q},\delta)$ is the polynomial ring $\Cx[s_{apq} \mid a \in \overline{Q},1 \leq p \leq \delta_{h(a)},1 \leq q \leq \delta_{t(a)} ]$ where the indeterminate $s_{apq}$ picks out the $p,q$ entry of the matrix $x_a$, corresponding to $x \in \uRep(\overline{Q},\delta)$. It is fairly straightforward to see that $\Vs_Q(\tau,\delta)$ has coordinate ring $R[s_{apq}]/J_{\tau}$, where $J_{\tau}$ is generated by the elements
$$
\sum_{h(a)=i} \sum_{r=1}^{\delta_{t(a)}}s_{apr}s_{a^*rq}-\sum_{t(a)=i} \sum_{r=1}^{\delta_{h(a)}}s_{a^*pr}s_{arq}-\delta_{pq}\tau_i
$$
for each vertex $i$ and for $1 \leq p,q \leq \delta_i$. Letting $\ell=\sum_i \delta_i$ there is a natural ring homomorphism $R \overline{Q} \fl \uMat_{\ell}(R[s_{apq}])$ sending an arrow $a$ to the matrix whose entries are the relevant $s_{apq}$. By our previous remark this homomorphism descends to a map $\mathfrak{A}^{\tau}(Q) \fl \uMat_{\ell}(\Cx[\Vs_Q(\tau,\delta)])$. Since $\delta_0=1$, this restricts to a homomorphism $e_0  \mathfrak{A}^{\tau}(Q) e_0 \fl \Cx[\Vs_Q(\tau,\delta)]$. One easily checks that the elements in the image of this map are invariant under the action of $\uGL(\delta)$. In this way we get a map $\phi_{\tau}\colon  C^{\tau} \fl \Cx[\Vs_Q(\tau,\delta)]^{\uGL(\delta)}$.  It follows from \cite[Corollary~8.12]{CrH98} that if $\tau \in R^I$ satisfies $\sum_{i}\delta_i \tau_i=0$, then the map $\phi_{\tau}$ is an isomorphism. Thus, we arrive to the following result.

\begin{prop}\label{prop3.4}
If $\sum_i \delta_i \tau_i=0$, then $C^{\tau} \cong \Cx[\Vs_Q(\tau,\delta)]^{\uGL(\delta)}$.
\end{prop} 

One immediate consequence of this is that $ \Cx[\Vs_Q(\tau,\delta)]^{\uGL(\delta)}$ is an integrally closed domain, so the quotient scheme $\Vs_Q(\tau,\delta)\sslash \uGL(\delta)$ is normal.

Another application of Proposition~\ref{prop3.4} is given by the following.

\begin{cor}
If $\sum_i \delta_i \tau_i=0$, then the rings $A^{\tau}$ and $C^{\tau}$ have Krull dimension $3$.
\end{cor}

\begin{proof}
Using Lemma~\ref{lem3.1} and \cite[Corollary~13.4.9]{McCR87} one sees immediately that the rings $A^{\tau}$ and $C^{\tau}$ are PI rings, and so their Krull dimension coincides with their GK dimension. The assertion follows.
\end{proof}

We finish this subsection with an observation which will be central to our main result. Here we denote the centres of $A^{\tau}$ and $C^{\tau}$ by $Z(A^{\tau})$ and $Z(C^{\tau})$ respectively.

\begin{prop}\label{prop3.6}
The map $\phi\colon  A^{\tau} \fl C^{\tau}$ given by $\phi(a)=eae$ for all $a$ in $A^{\tau}$ restricts to an algebra isomorphism from $Z(A^{\tau})$ to $Z(C^{\tau})$.
\end{prop}

\begin{proof}
It is a straightforward calculation to show that $\phi |_{Z(A^{\tau})}$ is an algebra homomorphism with image in $Z(C^{\tau})$. To see that it is an algebra isomorphism we construct the inverse map. First we note that an element $\xi$ in $Z(C^{\tau})$ implements a $C^{\tau}$-endomorphism of $A^{\tau}e$ via right multiplication by $\xi$. Thanks to Lemma~\ref{lem3.1}, this endomorphism can be regarded as an element $a_{\xi}$ of $A^{\tau}$. Then the algebra homomorphism $\psi \colon  Z(C^{\tau}) \fl A^{\tau}$ given by $\psi(\xi)=a_{\xi}$ for all $\xi$ in $Z(C^{\tau})$ has its image in $Z(A^{\tau})$ because the right multiplication by $\xi$ on $A^{\tau}e$ commutes with left multiplication by $A^{\tau}$. It is readily verified that this homomorphism is inverse to $\phi |_{Z(A^{\tau})}$. This completes the proof of the proposition. 
\end{proof}

This result shows a second vital feature of $C^{\tau}$:~its structure determines the centre of $A^{\tau}$. Now, if $\sum_i \delta_i \tau_i=0$, then we know from Proposition~\ref{prop3.4} that $C^{\tau}$ is commutative. According to Proposition~\ref{prop3.6}, in this case $C^{\tau} \cong Z(A^{\tau})$. 

\subsection{Brief account of Van den Bergh's construction}\label{subsec3.2}
In this subsection we describe some of Van den Bergh's results concerning noncommutative crepant resolutions.

We first consider the following more general situation. Let $R$ be a commutative noetherian algebra over $\Cx$ and let $A$ be an $R$-algebra which is finitely generated as an $R$-module. Let $\{e_0,\dots,e_n \}$ be a complete set of primitive orthogonal idempotents in $A$ and set $I=\{0,1,\dots,n \}$. We wish to construct a moduli space of $A$-modules. To do this we introduce a stability condition.

Let us fix a field $K$ and a ring homomorphism $R \fl K$. If $M$ is a finite dimensional $A \otimes_R K$-module, its dimension vector $\udimvec M$ is the element of $\Nx^I$ whose $i$th component is $\dim_K (e_iM)$. Let $\theta$ be a homomorphism $\Zx^I \fl \Zx$. A finite dimensional $A \otimes_R K$-module $M$ is said to be $\theta$-{\it stable} (or $\theta$-{\it semistable}) if $\theta(\udimvec M)=0$, but $\theta(\udimvec M')>0$ (or $\theta(\udimvec M')\geq 0$) for every proper submodule $M' \subseteq M$. As usual, we say that $\theta$ is {\it generic} for $\alpha$ if every $\theta$-semistable $A \otimes_R K$-module of dimension $\alpha$ is $\theta$-stable. Note that such a $\theta$ exists if and only if $\alpha$ is indivisible, meaning that the $\alpha_i$ have no common divisors. As a matter of fact, the condition $\theta(\beta) \neq 0$ for all $0 < \beta < \alpha$ ensures $\theta$ is generic. 

Next we recall the notion of family from \cite{GS04}. Fix a dimension vector $\alpha \in \Nx^I$. A family of $A$-modules of dimension $\alpha$ over an $R$-scheme $S$ is a locally free sheaf $\Fs$ over $S$ together with an $R$-algebra homomorphism $A \fl \uEnd_S(\Fs)$ such that $e_i \Fs$ has constant rank $\alpha_i$ for all $i \in I$. Two such families $\Fs$ and $\Fs'$ are equivalent if there is a line bundle $\Ls$ on $S$ and an isomorphism $\Fs \cong \Fs' \otimes_{\Os_S} \Ls$. Finally we say that a family $\Fs$ is $\theta$-{\it stable} (or $\theta$-{\it semistable}) if for every field $K$ and every morphism $\phi\colon  \uSpec K \fl S$ we have that $\phi^* \Fs$ is $\theta$-stable (or $\theta$-semistable) as $A \otimes_R K$-module. 

We have the following result, see \cite[Proposition~6.2.1]{Bergh02}.

\begin{prop}\label{prop3.8}
If $\theta$ is generic, then the functor which assigns to a scheme $S$ the set of equivalence classes of families of $\theta$-stable $A$-modules of dimension $\alpha$ over $S$ is representable by a projective scheme $\Ms_{\theta}(A,\alpha)$ over $X=\uSpec R$. 
\end{prop}

We now illustrate how to use this result to construct a crepant resolution starting from a noncommutative one. Let $R$ be an integrally closed Gorenstein domain admitting a noncommutative crepant resolution $A = \uEnd_R(M)$ and set $X=\uSpec R$. We assume, for simplicity, that $X$ is irreducible. Let $M= \bigoplus_{i \in I} M_i$ be any decomposition of $M$ corresponding to idempotents $e_0,\dots,e_n \in A = \uEnd_R(M)$, and let $\alpha \in \Nx^I$ be the vector with $\alpha_i = \urank M_i$. By Proposition~\ref{prop3.8} we know that, for generic $\theta$, there is a fine moduli space $\Ms_{\theta}(A,\alpha)$ of $\theta$-stable $A$-modules of dimension $\alpha$. Let us denote by $\phi\colon  \Ms_{\theta}(A,\alpha) \fl X$ the structure morphism. If we let $U \subseteq X$ be the open subset over which $M$ is locally free then it follows from \cite[Lemma~6.2.3]{Bergh02} that $\phi^{-1}(U) \fl U$ is an isomorphism. Each point $y \in \phi^{-1}(U)$ is a  $\theta$-stable $A$-module of dimension $\alpha$ so there is an embedding $U \hookrightarrow  \Ms_{\theta}(A,\alpha)$. Let $W \subseteq \Ms_{\theta}(A,\alpha)$ be the irreducible component of $\Ms_{\theta}(A,\alpha)$ containing the image of this morphism. Then $W$ is fine, in that $W$ is projective and there is a universal sheaf $\Us$ on $W \times X$. We denote by $\Ps$ the restriction of $\Us$ to $W$. Notice that $\Ps$ is a sheaf of (left) $A$-modules on $W$.

Now let $\Dd^b(\uCoh(W))$ denote the bounded derived category of coherent sheaves on $W$ and $\Dd^b(\umod A)$ the bounded derived category of finitely generated right modules over $A$. The method of Bridgeland, King and Reid generalises to prove the following result, see \cite[Theorem~6.3.1]{Bergh02}.

\begin{thm}\label{thm3.9}
Let the setting be as above. If $\dim (W \times_X W) \leq \dim X+1$, then $\phi \colon  W \fl X$ is a crepant resolution and the functors $\Rd \Gamma(- \otimes_{\Os_W}^{\Ld} \Ps)$ and $-\otimes_A^{\Ld}  \Rd \Hs om\spdot_{W}(\Ps,\Os_W)$ define inverse equivalences between $\Dd^b(\uCoh(W))$ and $\Dd^b(\umod A)$.
\end{thm}

\subsection{Application to our situation}
We now return to the concrete situation of Sect.~\ref{subsec3.1}. Our main aim is to show how the ideas developed in the previous subsection can be used to prove that any small resolution of an ADE fibered Calabi-Yau threefold is derived equivalent to the corresponding $N=1$ ADE quiver algebra. We start with some preliminary observations. 

We have seen in Proposition~\ref{prop3.2} that the algebra $A^{\tau}$ is a noncommutative crepant resolution of $C^{\tau}$. Hereafter we assume that $\tau \in R^I$ satisfies $\sum_{i}\delta_i \tau_i=0$. As we pointed out earlier, this implies that $C^{\tau} \cong Z(A^{\tau}) \cong \Cx[\Vs_Q(\tau,\delta)]^{\uGL(\delta)}$. Setting $X=\uSpec C^{\tau}$, it follows, from Theorem~\ref{thmNC3fold}, that $X$ is isomorphic to an ADE fibered Calabi-Yau threefold.

We shall again let $\rho_0,\dots,\rho_n$ denote the irreducible representations of $G$ with $\rho_0$ trivial, and set $I=\{ 0,\dots,n \}$. For each $i \in I$, let $f_i$ be the idempotent in $\Cx G$ with $\Cx G f_i \cong \rho_i$. As previously emphasized, we may regard the $f_i$'s as elements of $A^{\tau}$. Then $f_i A^{\tau}e$ is a submodule of $A^{\tau}e$ for all $i \in I$ and $A^{\tau}e=\bigoplus_{i \in I}f_i A^{\tau}e$. Bearing in mind that $f_i A^{\tau}e \cong \uHom_{\Cx G}(\rho_i, A^{\tau}e)$ we have $\delta_i=\dim \rho_i = \urank (f_i A^{\tau}e)$. Let $\Ms_{\theta}(A^{\tau},\delta)$ be the moduli space, as constructed in the previous subsection, of $\theta$-stable $A^{\tau}$-modules of dimension $\delta$ (equivalently, isomorphic to $\Cx G$), and let $W$ be the irreducible component of $\Ms_{\theta}(A^{\tau},\delta)$ that maps birationally to $X$. 

With the aid of Theorem~\ref{thm3.9}, we easily derive the following.

\begin{prop}\label{prop3.11}
With the notation above, $W$ is a crepant resolution of $X$ and there is an equivalence of categories between $\Dd^b(\uCoh(W))$ and $\Dd^b(\umod A^{\tau})$.
\end{prop}

\begin{proof}
Define $\Delta$ to be the diagonal of $W \times W$. As in the previous subsection, we write $\phi\colon  W \fl X$ for the structure morphism. This is a birational projective mapping, so it is closed. Let us take non-empty open subsets $V \subseteq W$ and $U \subseteq X$, such that $\phi$ restricts to an isomorphism $\phi\colon  V \fl U$. Denote by $Z$ the complement of $V$.  We may assume without loss of generality that $\phi(Z) \cap U = \varnothing$. It therefore follows that $W \times_X W \subseteq \Delta \cap (Z \times Z)$. Since $\dim X = 3$ we have $\dim Z \leq 2$, which ensures that $\dim (W \times_X W) \leq 4$. Now we are in the situation of Theorem~\ref{thm3.9} and the assertion follows. 
\end{proof}

We now apply Proposition~\ref{prop3.11} to prove the result promised in the beginning of this subsection.

\begin{thm}
Let the context be as above. If $\pi \colon  Y \fl X$ is a small resolution of the ADE fibered Calabi-Yau threefold $X$, then there is an equivalence of categories
$$
\Dd^{b}(\uCoh(Y)) \cong \Dd^b(\umod \mathfrak{A}^{\tau}(Q)),
$$
where $\mathfrak{A}^{\tau}(Q)$ is the associated $N=1$ ADE quiver algebra.
\end{thm}

\begin{proof}
It is well-known (see, e.g.~\cite[Proposition~16.4]{CKM88}) that $\pi$ is a crepant resolution. Owing to Proposition~\ref{prop3.11}, there exists another crepant resolution $\phi\colon  W \fl X$ associated to $A^{\tau}$. Let $f \colon  Y \fl W$ be the birational map over $X$ such that $f$ is isomorphic in codimension $1$. Then, by \cite[Theorem~6.38]{KoMo98}, $f$ is a composition of finitely many flops. A result of Bridgeland \cite[Theorem~1.1]{Bridg02} provides an equivalence of categories $\Dd^b(\uCoh(Y))\cong \Dd^b(\uCoh(W))$. Invoking Propositions \ref{prop3.11} and \ref{prop3.3}, we therefore deduce that
$$
\Dd^b(\uCoh(Y))\cong\Dd^b(\umod A^{\tau}) \cong \Dd^b(\umod \mathfrak{A}^{\tau}(Q)),
$$
as we wished to show. 
\end{proof}

\def\cprime{$'$} \def\cprime{$'$}

\addcontentsline{toc}{section}{Bibliography}

\begin{thebibliography}{10}

\bibitem{Asp05}
P.~S. Aspinwall.
\newblock D-branes on {C}alabi-{Y}au manifolds.
\newblock In {\em Progress in string theory}, pages 1--152. World Sci. Publ.,
  Hackensack, NJ, 2005.

\bibitem{Asp06}
P.~S. Aspinwall.
\newblock D-branes, {$\Pi$}-stability and {$\theta$}-stability.
\newblock In {\em Snowbird lectures on string geometry}, volume 401 of {\em
  Contemp. Math.}, pages 1--13. Amer. Math. Soc., Providence, RI, 2006.

\bibitem{BGLP00}
C.~Beasley, B.~R. Greene, C.~I. Lazaroiu, and M.~R. Plesser.
\newblock D3-branes on partial resolutions of abelian quotient singularities of
  {C}alabi-{Y}au threefolds.
\newblock {\em Nuclear Phys. B}, 566(3):599--641, 2000.

\bibitem{Ber02}
D.~Berenstein.
\newblock Reverse geometric engineering of singularities.
\newblock {\em J. High Energy Phys.}, (4):No. 52, 18, 2002.

\bibitem{Bergman08}
A.~Bergman.
\newblock Stability conditions and branes at singularities.
\newblock {\em J. High Energy Phys.}, (10):073, 19, 2008.

\bibitem{BBS03}
R.~Bocklandt, L.~Le~Bruyn, and S.~Symens.
\newblock Isolated singularities, smooth orders, and {A}uslander regularity.
\newblock {\em Comm. Algebra}, 31(12):6019--6036, 2003.

\bibitem{Bridg02}
T.~Bridgeland.
\newblock Flops and derived categories.
\newblock {\em Invent. Math.}, 147(3):613--632, 2002.

\bibitem{BKR01}
T.~Bridgeland, A.~King, and M.~Reid.
\newblock The {M}c{K}ay correspondence as an equivalence of derived categories.
\newblock {\em J. Amer. Math. Soc.}, 14(3):535--554 (electronic), 2001.

\bibitem{Briesk68}
E.~Brieskorn.
\newblock Die {A}ufl\"osung der rationalen {S}ingularit\"aten holomorpher
  {A}bbildungen.
\newblock {\em Math. Ann.}, 178:255--270, 1968.

\bibitem{CFIKV02}
F.~Cachazo, B.~Fiol, K.~Intriligator, S.~Katz, and C.~Vafa.
\newblock A geometric unification of dualities.
\newblock {\em Nuclear Phys. B}, 628(1-2):3--78, 2002.

\bibitem{CKV01}
F.~Cachazo, S.~Katz, and C.~Vafa.
\newblock Geometric transitions and {$N=1$} quiver theories.
\newblock http://arxiv.org/abs/hep-th/0108120.

\bibitem{CS98}
H.~Cassens and P.~Slodowy.
\newblock On {K}leinian singularities and quivers.
\newblock In {\em Singularities (Oberwolfach, 1996)}, volume 162 of {\em Progr.
  Math.}, pages 263--288. Birkh\"auser, Basel, 1998.

\bibitem{CKM88}
H.~Clemens, J.~Koll{\'a}r, and S.~Mori.
\newblock Higher-dimensional complex geometry.
\newblock {\em Ast\'erisque}, (166):144 pp. (1989), 1988.

\bibitem{CrIs04}
A.~Craw and A.~Ishii.
\newblock Flops of {$G$}-{H}ilb and equivalences of derived categories by
  variation of {GIT} quotient.
\newblock {\em Duke Math. J.}, 124(2):259--307, 2004.

\bibitem{Cr01}
W.~Crawley-Boevey.
\newblock Geometry of the moment map for representations of quivers.
\newblock {\em Compositio Math.}, 126(3):257--293, 2001.

\bibitem{CEG07}
W.~Crawley-Boevey, P.~Etingof, and V.~Ginzburg.
\newblock Noncommutative geometry and quiver algebras.
\newblock {\em Adv. Math.}, 209(1):274--336, 2007.

\bibitem{CrH98}
W.~Crawley-Boevey and M.~P. Holland.
\newblock Noncommutative deformations of {K}leinian singularities.
\newblock {\em Duke Math. J.}, 92(3):605--635, 1998.

\bibitem{DelFreed99}
P.~Deligne and D.~S. Freed.
\newblock Supersolutions.
\newblock In {\em Quantum fields and strings: a course for mathematicians,
  {V}ol. 1, 2 ({P}rinceton, {NJ}, 1996/1997)}, pages 227--355. Amer. Math.
  Soc., Providence, RI, 1999.

\bibitem{FHHZ08a}
D.~Forcella, A.~Hanany, Y.-H. He, and A.~Zaffaroni.
\newblock The master space of $\mathcal{N}=1$ gauge theories.
\newblock {\em J. High Energy Phys.}, (8):012, 72, 2008.

\bibitem{Gin06}
V.~Ginzburg.
\newblock Calabi-{Y}au algebras.
\newblock http://aps.arxiv.org/abs/math/0612139v2.

\bibitem{GS04}
I.~Gordon and S.~P. Smith.
\newblock Representations of symplectic reflection algebras and resolutions of
  deformations of symplectic quotient singularities.
\newblock {\em Math. Ann.}, 330(1):185--200, 2004.

\bibitem{KV00}
M.~Kapranov and E.~Vasserot.
\newblock Kleinian singularities, derived categories and {H}all algebras.
\newblock {\em Math. Ann.}, 316(3):565--576, 2000.

\bibitem{K04}
S.~Katz.
\newblock A{D}{E} geometry and dualities.
\newblock {\em Workshop on Algebraic Geometry and Physics, Lisbon}, 2004.

\bibitem{KM92}
S.~Katz and D.~R. Morrison.
\newblock Gorenstein threefold singularities with small resolutions via
  invariant theory for {W}eyl groups.
\newblock {\em J. Algebraic Geom.}, 1(3):449--530, 1992.

\bibitem{K94}
A.~D. King.
\newblock Moduli of representations of finite-dimensional algebras.
\newblock {\em Quart. J. Math. Oxford Ser. (2)}, 45(180):515--530, 1994.

\bibitem{KW99}
I.~R. Klebanov and E.~Witten.
\newblock Superconformal field theory on threebranes at a {C}alabi-{Y}au
  singularity.
\newblock {\em Nuclear Phys. B}, 536(1-2):199--218, 1999.

\bibitem{KoMo98}
J.~Koll{\'a}r and S.~Mori.
\newblock {\em Birational geometry of algebraic varieties}, volume 134 of {\em
  Cambridge Tracts in Mathematics}.
\newblock Cambridge University Press, Cambridge, 1998.
\newblock With the collaboration of C. H. Clemens and A. Corti, Translated from
  the 1998 Japanese original.

\bibitem{Kont93}
M.~Kontsevich.
\newblock Formal (non)commutative symplectic geometry.
\newblock In {\em The {G}el\cprime fand {M}athematical {S}eminars, 1990--1992},
  pages 173--187. Birkh\"auser Boston, Boston, MA, 1993.

\bibitem{Br99}
L.~Le~Bruyn.
\newblock Noncommutative compact manifolds constructed from quivers.
\newblock {\em AMA Algebra Montp. Announc.}, pages Paper 1, 5 pp. (electronic),
  1999.

\bibitem{Br06}
L.~Le~Bruyn.
\newblock Non-commutative algebraic geometry and commutative
  desingularizations.
\newblock In {\em Noncommutative algebra and geometry}, volume 243 of {\em
  Lect. Notes Pure Appl. Math.}, pages 203--252. Chapman \& Hall/CRC, Boca
  Raton, FL, 2006.

\bibitem{BrSym05}
L.~Le~Bruyn and S.~Symens.
\newblock Partial desingularizations arising from non-commutative algebras.
\newblock http://arxiv.org/abs/math.RA/0507494.

\bibitem{McCR87}
J.~C. McConnell and J.~C. Robson.
\newblock {\em Noncommutative {N}oetherian rings}.
\newblock Pure and Applied Mathematics (New York). John Wiley \& Sons Ltd.,
  Chichester, 1987.
\newblock With the cooperation of L. W. Small, A Wiley-Interscience
  Publication.

\bibitem{Seg08}
E.~Segal.
\newblock The {$A\sb \infty$} deformation theory of a point and the derived
  categories of local {C}alabi-{Y}aus.
\newblock {\em J. Algebra}, 320(8):3232--3268, 2008.

\bibitem{Slodow80}
P.~Slodowy.
\newblock {\em Simple singularities and simple algebraic groups}, volume 815 of
  {\em Lecture Notes in Mathematics}.
\newblock Springer, Berlin, 1980.

\bibitem{StaZha94}
J.~T. Stafford and J.~J. Zhang.
\newblock Homological properties of (graded) {N}oetherian {${\rm PI}$} rings.
\newblock {\em J. Algebra}, 168(3):988--1026, 1994.

\bibitem{Sz04}
B.~Szendr{\H{o}}i.
\newblock Artin group actions on derived categories of threefolds.
\newblock {\em J. Reine Angew. Math.}, 572:139--166, 2004.

\bibitem{Sz07}
B.~Szendr{\H{o}}i.
\newblock Non-commutative {D}onaldson-{T}homas invariants and the conifold.
\newblock {\em Geom. Topol.}, 12(2):1171--1202, 2008.

\bibitem{Sz08}
B.~Szendr{\H{o}}i.
\newblock Sheaves on fibered threefolds and quiver sheaves.
\newblock {\em Comm. Math. Phys.}, 278(3):627--641, 2008.

\bibitem{Tj70}
G.~N. Tjurina.
\newblock Resolution of singularities of flat deformations of double rational
  points.
\newblock {\em Funkcional. Anal. i Prilo\v zen.}, 4(1):77--83, 1970.

\bibitem{Bergh02}
M.~van~den Bergh.
\newblock Non-commutative crepant resolutions.
\newblock In {\em The legacy of Niels Henrik Abel}, pages 749--770. Springer,
  Berlin, 2004.

\bibitem{Bergh04}
M.~Van~den Bergh.
\newblock Three-dimensional flops and noncommutative rings.
\newblock {\em Duke Math. J.}, 122(3):423--455, 2004.

\bibitem{Wem07}
M.~Wemyss.
\newblock Reconstruction algebras of type {A}.
\newblock http://arxiv.org/abs/0704.3693.

\bibitem{Wijn08}
M.~Wijnholt.
\newblock Parameter space of quiver gauge theories.
\newblock {\em Adv. Theor. Math. Phys.}, 12(4):711--755, 2008.

\bibitem{Zhu06}
X.~Zhu.
\newblock Representations of {$N=1$} {$ADE$} quivers via reflection functors.
\newblock {\em Michigan Math. J.}, 54(3):671--686, 2006.

\end{thebibliography}
\end{document}